\documentclass{amsart}
\usepackage{amsfonts}
\usepackage{amsthm}
\usepackage[usenames]{color}
\usepackage{amssymb}
\usepackage{amsmath}
\usepackage{amscd}

 \usepackage{paralist}
  \usepackage[colorlinks=true]{hyperref}
    \hypersetup{urlcolor=blue, citecolor=red}

\newtheorem{thm}{Theorem}[section]
\newtheorem{lem}[thm]{Lemma}
\newtheorem{prop}[thm]{Proposition}
\newtheorem{cor}[thm]{Corollary}

\theoremstyle{definition} \newtheorem*{defn}{Definition}
\theoremstyle{definition} \newtheorem*{rem}{Remark}

\newcommand{\mb}[1]{\mathbb{#1}}
\newcommand{\mc}[1]{\mathcal{#1}}

\title[Comparison for Sobolev gradient flow]{A comparison principle
for a Sobolev gradient semi-flow}
\author[T. Blass, R. de la Llave and E. Valdinoci]{}
\subjclass{Primary: 35B50, 46N20; Secondary: 35J20}
\keywords{Comparison principle, Sobolev gradient, semigroups of linear
operators, fractional powers of elliptic operators}

\email{tblass@math.utexas.edu}
\email{llave@math.utexas.edu}
\email{enrico@math.utexas.edu}
\begin{document}

\maketitle
%=====================
\centerline{\scshape Timothy Blass and Rafael de la Llave}
\medskip
{\footnotesize
% please put the address of the first author
 \centerline{Department of Mathematics, 1 University Station C1200}
   \centerline{Austin, TX 78712-0257, USA}
  % \centerline{ Springfield, MO 65801-2604, USA}
} % Do not forget to end the {\footnotesize by the sign }

\medskip

\centerline{\scshape Enrico Valdinoci}
\medskip
{\footnotesize
 % please put the address of the second  and third author
\centerline{Dipartimento di Matematica, Universit\`a di Roma Tor Vergata}
   \centerline{Via della Ricerca Scientifica, 1, I-00133 Roma, Italy}
}
\bigskip

% The name of the associate editor will be entered by an editorial staff
% \centerline{(Communicated by Xavier Cabr\'e)}

%======================
\begin{abstract}
We consider gradient descent equations for energy functionals of
the type $S(u) = \frac{1}{2}\langle u(x), A(x)u(x) \rangle_{L^2} +
\int_{\Omega} V(x,u) \, dx$, where $A$ is a
uniformly elliptic operator of order 2, with smooth coefficients.
The gradient descent equation for such a functional depends on the
metric under consideration.

We consider the steepest descent equation
for $S$  where the gradient is an element of the Sobolev space
$H^{\beta}$, $\beta \in (0,1)$, with a metric that depends on $A$
and a positive number $\gamma > \sup |V_{22}|$.
We prove a weak comparison principle for such a gradient flow.

We extend our methods to the case where $A$ is a fractional power of
an elliptic operator, and
provide an application to the Aubry-Mather theory for partial
differential equations and pseudo-differential equations by
finding plane-like minimizers of the energy functional.

\end{abstract}
%\subjclass{}
%\footnote{
%{\bf Subject Classification:}
%35B50,***  %Maximum principles
%46N20,***  %Applications to differential and integral equations
%35A15,  %Variational methods
%35J20***  %Variational methods for second-order, elliptic equations
%35B27,  %Homogenization; partial differential equations in media with periodic structure
%45P05,  %Integral operators
%58E15,  %Application to extremal problems in several variables
%47D03  %Groups and semigroups of linear operators
%}
\section{Introduction}
In this paper we prove a comparison principle for steepest descent equations,
in the Sobolev gradient direction (see (\ref{eqn:general})). When one is
interested in minimizing functionals of the type
\begin{equation} \label{eqn:enfunc}
S(u) = \frac{1}{2}\langle u(x), A(x)u(x) \rangle_{L^2} +
\int_{\Omega} V(x,u) \, dx,
\end{equation}
where $A$ is an elliptic operator,
it is natural to consider the gradient descent equation $\partial_tu =
-\nabla S(u)$.
The gradient of $S$ depends on the metric under consideration.
Our main result is a comparison principle for flows of this type, which
is formulated in Section (\ref{sec:main}), where $\nabla S(u)$ is an
element of the Sobolev space $H^{\beta}$, $\beta \in (0,1)$.
The methods used to prove the
comparison principle may be of independent interest and are outlined at
the end of Section
(\ref{sec:main}). In particular, the methods extend naturally to
show a comparison result for $\partial_tu = -\nabla S(u)$, with $A$
replaced  by a fractional power $A^{\alpha}$, $\alpha \in (0,1)$, as
shown in Section \ref{sec:feqns}.  The metrics
and Sobolev gradients we consider are explained in Sections \ref{sec:space}
and \ref{sec:grad}. See also \cite{Neuberger97}.

More concretely, we consider a self-adjoint, uniformly elliptic operator $A$
given by the formula
\begin{equation} \label{eqn:A}
Au = - \sum_{i,j=1}^d \partial_{x_j}\left(a^{ij}(x)\partial_{x_i}u\right) =
-\mbox{div} (a(x) \nabla u),
\end{equation}
where the coefficient functions, $a^{ij} \in C^{\infty}(\mb{R}^d)$ are
symmetric in $i,j$ and we have positive constants  $\Lambda_1,\, \Lambda_2$
such that for every$x \in \mb{R}^d$
\[
\Lambda_1 |\xi|^2 \leq \sum_{i,j=1}^d a^{ij}(x)\xi_i \xi_j \leq
\Lambda_2 |\xi|^2, \quad \forall \xi \in \mb{R}^d.
\]
Then, for a suitably large constant $\gamma$, and for $\beta \in (0,1)$,
we will show a comparison principle for the flow defined by the evolution
equation
\begin{equation} \label{eqn:general}
\partial_t u = -(\gamma + A)^{1-\beta}u + (\gamma +A)^{-\beta}\left(\gamma u -
V_2(x,u) \right),
\end{equation}
where $V(x,y)\in C^r(\mb{R}^d \times \mb{R})$, $r \geq 2$, and $V_2$ denotes
the derivative of $V$ with respect to its last argument. The fractional
powers of $\gamma + A$ that appear in (\ref{eqn:general}) will be defined in
Section \ref{sec:pow}. As we shall show over the next few sections of this
introduction, equation (\ref{eqn:general}) is the steepest descent equation
for $S$ in the Sobolev space $H^{\beta}$ with inner product
$\langle u,v \rangle_{H^{\beta}}=\langle (\gamma + A)^{\beta} u,v \rangle_{L^2}$,
as explained in Section \ref{sec:space}. The domain and boundary considerations
for equation (\ref{eqn:general}) are discussed in the following section.
A sufficient lower bound for the constant $\gamma$ will be given in Section
\ref{sec:comp} and will depend on the nonlinear term $V$.

In Section
\ref{sec:main} we state the comparison result in its full generality but
leave the proof until Section \ref{sec:comp}. We gather some
previous results in Sections \ref{sec:Moser} and \ref{sec:pow}
and apply them to our problem in Section \ref{sec:appl}. Section
\ref{sec:eu} is devoted to the proofs of existence
and uniqueness of solutions to equation (\ref{eqn:general}).
If $S$ were $C^2$, then a theorem in \cite{Neuberger97}
would give existence and uniqueness immediately. However,
the functional $S$ defined in (\ref{eqn:enfunc}) is not even continuous
from $H^{\beta}$ to $\mb{R}$ for $\beta<1$.

As mentioned earlier, Section \ref{sec:feqns}
is devoted to explaining how
the techniques developed in the proof
of the comparison principle can be applied to fractional powers of
elliptic operators as well. That is, the techniques apply to equations
with the same form as (\ref{eqn:general}) with the operator $A$ replaced
by $A^{\alpha}$, for $\alpha \in (0,1)$.

Finally, we present an application to Aubry-Mather
theory for PDEs and pseudo-DEs in Section \ref{sec:app}. Aubry-Mather
theory concerns the minimizers of Lagrangian actions,
which can be classified by their associated rotation
vectors (or frequency vectors). In the PDE setting, Moser
extended the theory to certain types of energy functionals,
including those of the form (\ref{eqn:enfunc}), see \cite{Moser86}.
A certain geometric property (called the Birkhoff property)
of the minimizers is important to showing
existence of solutions for any rotation vector, see Section (\ref{sec:app}).
The use of gradient descent in this setting was introduced in
\cite{delaLlaveVal08}, where a comparison principle for the flow is
crucial to showing the solutions have this Birkhoff property.

%====================================================================
\subsection{Boundary conditions}
Our main results apply to Dirichlet boundary conditions for domains
$\Omega \subset \mb{R}^d$ that are compact with smooth boundary, as well as
for periodic boundary conditions (i.e.
$\Omega = \mb{T}^d \cong \mb{R}^d/\mb{Z}^d$).

In the
application in Section \ref{sec:app} we will also consider $\Omega = N\mb{T}^d
\cong \mb{R}^d/N \mb{Z}^d$, for which the reasoning regarding $N=1$ applies.
More succinctly, we pose (\ref{eqn:general}) as an initial-boundary value
problem, with $u(0,x) = u_0(x)\in L^{\infty}(\Omega)$ and one of the following
two cases:
\begin{eqnarray}
 u=0 &  &\mbox{on} \partial \Omega \label{eqn:dir} \\
 u(x+e,t) = u(x,t) & & \forall e \in \mb{Z}^d \label{eqn:per}
\end{eqnarray}
In the periodic setting we will require the functions $a^{ij}(x)$ and $V(x,y)$
to have period one in all variables.

Since our results rely on arguments that apply to both the Dirichlet and
periodic settings, we will not distinguish between the different cases
when stating the results.

%====================================================================
\subsection{Sobolev spaces} \label{sec:space}
There are several equivalent definitions of the Sobolev spaces $H^s(\Omega)$,
$s \in \mb{R}$. For the integer case,
we take $H^m(\Omega) = \{ u \in L^2(\Omega) : D^{\alpha}u \in L^2(\Omega), \,
\forall |\alpha| \leq m \} $, then the
intermediate spaces may be defined by interpolation methods \cite{Lunardi95},
\cite{Taylor96}. Alternatively, one can use the Fourier transform to define
$H^s(\mb{R}^d)= \{ u \in L^2(\mb{R}^d) : (1+|\xi|^2)^{s/2}\hat{u} \in
L^2(\mb{R}^d) \}$ and then $H^s(\Omega)=
\{ u \lvert_{\Omega} : u \in H^s(\mb{R}^d) \}$. The case where
$\Omega = \mb{T}^d$ is handled simply by replacing the Fourier transform
with the Fourier series. The factor $(1+|\xi|^2)^{s/2}$ makes it clear
that the operator $I-\Delta$ is the foundation of these spaces. Indeed,
$u\in H^s(\mb{R}^d)$ if and only if $(I-\Delta)^{s/2}u \in L^2(\mb{R}^d)$,
where $(I-\Delta)^{-s/2}$ is a particular case of the general definition
of the power of an elliptic operator given in equation (\ref{eqn:fpow1}),
and $(I-\Delta)^{s/2}$ is the inverse of $(I-\Delta)^{-s/2}$.

In fact, because
$\gamma + A$ is an (order 2) elliptic, self-adjoint operator, we can define
the Sobolev space $H_{\gamma, A}^s$ in the same manner as above, but replacing
$I-\Delta$ with $\gamma + A$, for $\gamma > 0$.
The inner product on $H_{\gamma, A}^s$ is given by
$\langle u,v \rangle_{s,\gamma,A} = \langle (\gamma + A)^s u,v \rangle_0 \equiv
\langle (\gamma + A)^s u,v \rangle_{L^2}$. In \cite{Shubin78}, page 57 it is
shown that the topology on $H_{\gamma, A}^s$ generated by the norm
obtained from the above inner
product is identical to the standard topology on $H^s$ (i.e. the topology
generated by
$\langle u,v \rangle_{H^s} \equiv \langle (I - \Delta)^s u,v \rangle_{L^2} $).
Thus, we henceforth omit the subscripts $\gamma$ and $A$ when referring to
$H_{\gamma, A}^s$ and $\langle \cdot ,\cdot \rangle_{s,\gamma,A}$.
Will will write $\langle \cdot ,\cdot \rangle_s$ for the inner product on $H^s$
and $\| \cdot \|_s$ for the norm on $H^s$. Though the topologies are
equivalent, the gradient of $S$ depends on the chosen inner product.
Thus, the gradient flow and therefore the comparison principle depend
on the choice of inner product.

%===========================================================
\subsection{Fractional powers of elliptic operators}
For $s > 0$, the operator $(\gamma + A)^{-s}$ is self-adjoint, bounded, linear,
invertible from $H^r$ to $H^{r+2s}$, and is defined by
\begin{equation} \label{eqn:fpow1}
(\gamma + A)^{-s} = \frac{1}{2 \pi i} \int_{\Gamma} z^{-s}
(\gamma + A-z)^{-1}\, dz,
\end{equation}
where $\Gamma$ is a rectifiable curve in the resolvent set $\rho(\gamma + A)
\subset \mb{C}$, avoiding $(-\infty,0]$. Here
$z^s$ is taken to be positive for positive real values of $z$ (see
\cite{Pazy83} page 69,
\cite{Shubin78} pages 83, 94). Positive powers are defined as
$(\gamma + A)^s = (\gamma + A)^k(\gamma + A)^{s-k}$ where $k \in \mb{N}$ and
$s < k$. It can be shown that $(\gamma + A)^s (\gamma + A)^r =
(\gamma + A)^{s+r}$ for $s$,$r \in \mb{R}$, and that if $s\in \mb{Z}$ our
definition coincides with the usual definition of integer powers of
$(\gamma + A)$, see \cite{Shubin78}. In particular we have
$(\gamma + A)^s = ((\gamma+A)^{-s})^{-1}$.

We will be interested in powers $0<\alpha < 1$, for which the integral in
(\ref{eqn:fpow1}) is equivalent to
\begin{equation} \label{eqn:fpow2}
(\gamma + A)^{-\alpha} = \frac{\sin{\pi \alpha}}{\pi} \int_0^{\infty} t^{-\alpha}
(t + \gamma + A)^{-1}\, dt.
\end{equation}
as shown in \cite{Pazy83}, Section 2.6. This fact will be needed in Section
\ref{sec:pow} when we discuss the semigroup theory related to
$(\gamma +A)^{1-\beta}$ and $(\gamma +A)^{-\beta}$.

We will use repeatedly in
Section \ref{sec:eu} that the operator $(\gamma +A)^{-\beta} \in
\mathcal{L}(H^s,H^{s+2\beta})$ is smoothing. We denote by
$\mathcal{L}(H_1,H_2)$
the space of bounded linear operators from the Hilbert space $H_1$ to the
Hilbert space $H_2$. For notational convenience we will
sometimes use $\lambda$ in place of $1-\beta$, in particular when describing
the domain $H^{s+2\lambda}$ of $(\gamma +A)^{\lambda}=(\gamma +A)^{1-\beta}$.

%====================================================================
\subsection{The Sobolev gradient} \label{sec:grad}
The motivation for equation (\ref{eqn:general}) is the desire to solve the
semilinear elliptic equation
\begin{equation} \label{eqn:ell}
-Au = V_2(x,u) \quad x\in \Omega
\end{equation}
subject to one of the boundary conditions (\ref{eqn:dir}) or (\ref{eqn:per}).
For background on equations of this type see
\cite{GilbargTrudinger01} and \cite{LadyUral68}.
Equation (\ref{eqn:ell}) has an associated variational principle. In fact, it
is the Euler-Lagrange equation for the functional
\begin{equation} \label{eqn:vp}
 S(u) = \int_{\Omega}\frac{1}{2} \left( a(x) \nabla u(x) \cdot \nabla u(x)
\right) + V(x,u(x))\, dx.
\end{equation}
To minimize $S$, and therefore find a solution to (\ref{eqn:ell}), we could
consider the steepest descent equation
\begin{equation} \label{eqn:sd}
\partial_t u = -Au - V_2(x,u).
\end{equation}
The motivation for equation (\ref{eqn:sd}) is that the right-hand side,
changed of sign, is
the unique element $g \in L^2$ such that $DS(u)\eta = \langle g,\eta
\rangle_{L^2}$, where $DS(u)$ is the Fr\'{e}chet derivative of $S$ at $u$.
This element $g \in L^2$ is called the $L^2$ gradient of $S$ with respect
to the inner product $\langle \cdot , \cdot \rangle_{L^2}$. Instead of
the standard $L^2$ inner product, if we used a different inner product, we
would obtain a different gradient for $S$.

We consider the Sobolev space $H^{\beta}$ with inner product
$\langle u,v \rangle_{\beta} = \langle (\gamma +A)^{\beta} u,v \rangle_{L^2}$,
and look for the gradient of $S$ with respect to this space and inner
product. That is, the unique element $g \in H^{\beta}$ such
that $DS(u)\eta = \langle g,\eta \rangle_{\beta}$.
We refer to $g$ as the Sobolev gradient of $S$ and write
$g=\nabla_{\beta}S(u)$ (see \cite{Neuberger97}). As noted at the end of
Section \ref{sec:space}, this gradient depends not only on $\beta$ but also
on our choice of inner product, which was determined by $A$ and $\gamma$.

We note that in each case, Dirichlet or periodic boundary conditions, we are
able to use the integration by parts formula
\begin{equation}
 -\int_{\Omega} \mbox{div} (a(x) \nabla u) v \, dx = \int_{\Omega} a(x) \nabla u
\cdot \nabla v \, dx,
\end{equation}
and we calculate the $H^{\beta}$-gradient as follows:
\begin{eqnarray*}
DS(u)\eta &=& \int_{\Omega} a(x)\nabla u\cdot \nabla \eta + V_2(x,u)\eta \,dx\\
&=&  \left\langle -\mbox{div}(a(x) \nabla u)+V_2(x,u), \eta \right\rangle_{L^2}
= \left\langle Au + V_2(x,u), \eta \right\rangle_{L^2} \\
&=& \left\langle (\gamma+A)^{\beta}(\gamma+A)^{-\beta}(A u + V_2(x,u)) ,\eta
\right\rangle_{L^2} \\
&=&  \left\langle (\gamma+A)^{-\beta}(\gamma u +A u -\gamma u + V_2(x,u)) ,\eta
\right\rangle_{{\beta}} \\
&=&  \left\langle (\gamma+A)^{1-\beta}u - (\gamma+A)^{-\beta}
(\gamma u - V_2(x,u)) , \eta \right\rangle_{{\beta}}.
\end{eqnarray*}
Thus, our steepest descent equation in $H^{\beta}$, $\partial_t u =
-\nabla_{\beta} S(u)$, becomes
\[
\partial_t u = -(\gamma+A)^{1-\beta}u + (\gamma+A)^{-\beta}(\gamma u - V_2(x,u)),
\]
which is identical to (\ref{eqn:general}). If the solution $u(x,t)$ of
(\ref{eqn:general}) approaches a critical point, that is $u(x,t) \to u_c(x)$ as
$t \to \infty$, then $u_c$ will solve $(\gamma+A)^{1-\beta}u_c =
(\gamma+A)^{-\beta}(\gamma u_c - V_2(x,u_c))$, which reduces to (\ref{eqn:ell}).

%====================================================================
\section{Main result} \label{sec:main}
We now wish to formulate our main theorem, which is a comparison principle
for the flow defined by (\ref{eqn:general}). The theorem is actually two
theorems, one for each type of boundary condition. Thus, in the statement of
the theorem, the space $H^s$ may refer to either of the two types of
Sobolev spaces $H^s_0$ (Dirichlet boundary conditions), or $H^s_P$
(periodic boundary conditions). We will write $\Omega$ to represent the
space domain, whether it is $\mb{T}^d$ or a smooth, bounded subset of
$\mb{R}^d$ with Dirichlet boundary conditions.
\begin{thm} \label{thm:1}
Let $V \in C^{r}(\Omega \times \mb{R},\mb{R})$, $r\geq 2$, and choose $\gamma$
such that $\gamma > \sup_{(x,y)\in \mb{R}^d \times \mb{R}} |V_{22}|$. Let $T>0$,
and let $u(x,t)$ and $v(x,t)$ be solutions of (\ref{eqn:general}) for
$t \in [0,T]$ with initial conditions $u(x,0) =u_0 \in L^{\infty}$ and
$v(x,0)=v_0 \in L^{\infty}$, where the exponent $\beta$ in (\ref{eqn:general})
is taken in the range $\beta \in (0,1)$. If $u_0 \geq v_0$ for almost every
$x\in \Omega$, then  $u(x,t) \geq v(x,t) \,\,\,
\forall t \in [0,T]$ and almost every $x \in \Omega$.
\end{thm}
\begin{rem}
The regularity of $V$ will be the only limit for the regularity of the
solution $u$. We show in Proposition \ref{prop:ext} that
$u(t,\cdot)\in H^{r+\delta-1}$, for any $\delta <2$. %Ref wants r+\delta-1 here
In particular, if $r > d/2-1$ then for
$u_0 \in L^{\infty}$ we have that the solution $u(t,\cdot)\in C^0(\Omega)$.
The proof of existence will
show that, in our case, we have $T=\infty$. This allows this semi-flow to be
used for finding critical points of $S$.
\end{rem}
We use the notation $L := -(\gamma + A)^{1-\beta}$ and
$X(u):=(\gamma + A)^{-\beta}(\gamma u - V_2(x,u))$ so that we can rewrite
(\ref{eqn:general}) as
\begin{equation} \label{eqn:simple}
\partial_t u = Lu + X(u),
\end{equation}
and easily refer to the linear and nonlinear components of the equation as
$L$ and $X$.
We now briefly outline the strategy for the proof of Theorem \ref{thm:1}.

We will show that $L$ generates a semigroup and this semigroup satisfies
a comparison principle. The theory of semigroups will also allow us to show
that, for large enough $\gamma$, the nonlinear operator $X$ will also satisfy
a comparison principle.

We will then show that solutions to equation
(\ref{eqn:simple}) exist for all time, and can be expressed via the
integral formula
\begin{equation} \label{eqn:duhamel}
u(x,t) = e^{tL}u_0(x) + \int_0^t e^{(t-s)L}X(u(s,x)) \, ds,
\end{equation}
commonly referred to as Duhamel's formula (see \cite{Taylor97}, page 272).
This is done by first proving $L^{\infty}$ estimates
for $e^{tL}$ and $X$ in Section \ref{sec:bounds}, which follow from
the comparison principles for each operator, respectively.
Then, using some results from Section
\ref{sec:appl}, we show that if $u$ belongs to a Sobolev space $H^{\sigma}$,
with $0 \leq \sigma <r+2$, then $e^{(t-s)L}X(u(s,x))$ belongs
to a higher space $H^{\sigma+\tau}$ with $\tau >0$ ($r$ is the regularity of
$V_2$).
The representation
in (\ref{eqn:duhamel}) then allows us to show $u \in H^{\sigma+\tau}$. This
is carried out in detail in Section \ref{sec:eu}.

To prove the comparison principle for solutions to
(\ref{eqn:simple}) we build and iteration scheme around
formula (\ref{eqn:duhamel}), namely:
$u^{j+1}(x,t) = e^{tL}u_0(x) + \int_0^t e^{(t-s)L}X(u^j(s,x)) ds$.
The comparison principles for $e^{tL}$ and $X$ will allow us
to show that $u_0 \geq v_0$ implies $u^{j}(x,t) \geq v^{j}(x,t)$.
Then we show that the $u^j$ converge to a solution of (\ref{eqn:simple})
that must also satisfy a comparison principle. This is done in
Section \ref{sec:comp}.

%====================================================================
\section{Preliminaries} \label{sec:prelim}
In Sections \ref{sec:Moser} and \ref{sec:pow} we
present previous results that will be applied in Section \ref{sec:appl}
to produce several results, including the comparison principles for $e^{tL}$
and $X$. In Section \ref{sec:bounds} we use the comparison results to produce
$L^{\infty}$ bounds on $e^{tL}$ and $X$, which will be important in proving
existence, uniqueness, and the final comparison result.
In particular, they allow the application
of the Moser estimates (\ref{eqn:me1}) and (\ref{eqn:me2}) below.

\subsection{Moser estimates} \label{sec:Moser}
The composition $V(x,u)$ will be controlled by the use of Moser estimates
for the composition of functions in $H^s$, $s \in \mb{N}$. If
$f\in C^s(\mb{R}^d\times \mb{R})$ and $\phi \in H^s(\mb{R}^d) \cap L^{\infty}
(\mb{R}^d)$ then
\begin{equation} \label{eqn:me1}
\|f(x,\phi)\|_s \leq c_s|f|_{C^s}(1+\|\phi \|_s).
\end{equation}
We also have that if $f\in C^{s+1}$ and if $\phi, \psi \in H^s$ are bounded
with $s\in \mb{N}$, then
\begin{equation} \label{eqn:me2}
\| f(x,\phi)-f(x,\psi) \|_s \leq c_s |f|_{C^{s+1}}
(1+\| \phi \|_s + \| \psi \|_s)\| \phi - \psi \|_s.
\end{equation}
The constant $c_s$ depends on the supremum of $\phi$ and the diameter of
$\Omega$, see \cite{Moser88},
\cite{Moser66}. When proving the existence of solutions to equation
(\ref{eqn:general}) we will first show that they exist in $L^{\infty}$ for
all time $t>0$, and then that they are in fact continuous in the domain
$\Omega$ and differentiable in time.

%================================================
\subsection{Properties of semigroups and fractional powers} \label{sec:pow}
In this section we gather some general bounds and properties
of semigroups generated by a class of operators called \emph{m-accretive}, and
in some cases self-adjoint m-accretive operators. Though these results are not
new,
we include them here because they are very useful and will be applied to
$L=-(\gamma +A)^{1-\beta}$ in Section \ref{sec:appl}.
\begin{defn}
If $H$ is a Hilbert space and $D \subset H$ is a dense linear subspace of $H$,
and if a linear operator $B: D\to H$ satisfies
\begin{equation} \label{eqn:accretive}
\begin{split}
 &\langle -Bu,u \rangle \geq  0 \quad \forall u\in D, \,\,\, \mbox{and} \\
&(-B+I)D =  H, \qquad \qquad \quad
\end{split}\end{equation}
then $-B$ is called \emph{m-accretive} (and $B$ is called
\emph{m-dissipative}).
\end{defn}
The Lumer-Phillips theorem (see \cite{Pazy83}) states that if $-B$
is m-accretive, then $B$
generates a strongly continuous semigroup of contractions, $e^{tB}$. That is,
$e^{tB} \in C([0,\infty),H)\cap C^1((0,\infty),H)$, there exists $c\geq 0$ such
that $\|e^{tB}\|_{\mathcal{L}(H)} \leq e^{-c t}$, and if $u_0 \in H$, then
$u(t,x):=e^{tB}u_0(x)$ satisfies
\begin{eqnarray*}
 \frac{\partial u}{\partial t} &=& Bu \\
u(0,x) &=& u_0(x).
\end{eqnarray*}

If $c>0$, then the fractional power $(-B)^{-\alpha}$ for $\alpha \in (0,1)$ as
defined in (\ref{eqn:fpow2}) can be expressed in terms of the semigroup
$e^{tB}$. We have the formula
\begin{equation} \label{eqn:fpowB}
 (-B)^{-\alpha}f = \frac{1}{\Gamma (\alpha)}\int_0^{\infty}t^{\alpha -1}
e^{tB}f \, dt,
\end{equation}
as shown in \cite{Pazy83}, \cite{Vrabie03}.

If $H$ is a Hilbert space and $B$ is self-adjoint and m-accretive on $H$
(this implies $B$ is regular m-accretive, as defined on page 22 of
\cite{Showalter97}), then for any integer
$n\geq 1$ and any $u \in H$, we have $e^{tB}u \in D(B^n)$ and that
\begin{equation} \label{eqn:B1}
\|B^n e^{tB}\|_{\mathcal{L}(H)} \leq \left( \frac{n}{\sqrt{2}t}\right)^n,
\end{equation}
see page 29 of \cite{Showalter97}.
This result applies to fractional powers of $-B$. In fact, if $B$ is as
above, and $\alpha \in (0,1)$ then there exists a
constant $C_{\alpha,T}$ such that for any $t\in (0,T]$ one has
\begin{equation}\label{eqn:B2}
\| (-B)^{\alpha} e^{tB} \|_{\mathcal{L}(H)} \leq C_{\alpha,T}
\left( \frac{\alpha}{t}\right)^{\alpha}.
\end{equation}
Furthermore, if we set
$Y=D((-B)^{\alpha})$, the domain of $(-B)^{\alpha}$, endowed with the graph
norm $\|u\|_Y=\|u\|_H +\|(-B)^{\alpha}u\|_H$, then
\begin{equation}\label{eqn:B3} \begin{split}
&\| e^{tB} \|_{\mathcal{L}(H,Y)} \leq
C_{\alpha,T} \left( \frac{\alpha}{t}\right)^{\alpha},\,\,\, \mbox{and}\\
&\| e^{tB}-I \|_{\mathcal{L}(Y,H)} \leq C'_{\alpha,T}t^{\alpha}.
\end{split} \end{equation}
For further details on (\ref{eqn:B1}), (\ref{eqn:B2}), and (\ref{eqn:B3})
see Section 4.1 of \cite{delaLlaveVal08}.

Finally, we will use the subordination identity of Bochner
\cite{Bochner49}. For $\sigma >0$, $t>0$, $\tau \geq 0$, and $0<\alpha <1$
we define
\begin{equation} \label{eqn:phi}
\phi_{t,\alpha}(\tau) = \frac{1}{2\pi i}\int_{\sigma - i\infty}^{\sigma + i\infty}
e^{\tau z - tz^{\alpha}}dz
\end{equation}
and $\phi_{t,\alpha}(\tau)=0$ if $\tau <0$. As is shown in \cite{Yosida74}
Section IX.11,
$\phi_{t,\alpha}(\tau)\geq 0$ for all $\tau >0$, and if $-B$ m-accretive,
we can represent $e^{-t(-B)^{\alpha}}$ as
\begin{equation} \label{eqn:BochGen}
e^{-t(-B)^{\alpha}}  =  \int_0^{\infty} e^{\tau B}
\phi_{t,\alpha}(\tau) \, d\tau, \quad t>0.
\end{equation}

%=====================================================
\subsection{Representation and comparison for $e^{tL}$ and $X$} \label{sec:appl}

We will show that the operator $L = -(\gamma + A)^{1-\beta}$ generates a
semigroup and that this semigroup has many nice properties, including a
comparison principle. We also show that the nonlinear operator $X$
satisfies a comparison principle.
Most of these facts will be derived from properties of
the semigroup generated by $-(\gamma + A)$.

\begin{prop}
For each $s\geq 0$, the operator $-(\gamma+A)$ is m-accretive with
respect to
$H^s$, and therefore generates a semigroup $e^{-(\gamma+A)t} \in
C([0,\infty),H^s)\cap C^1((0,\infty),H^s)$.
Moreover, the fractional power $(\gamma + A)^{-\alpha}$ for
$\alpha \in (0,1)$ as
defined in (\ref{eqn:fpow2}) can be expressed as
\begin{equation} \label{eqn:fpow3}
 (\gamma + A)^{-\alpha}f = \frac{1}{\Gamma (\alpha)}\int_0^{\infty}t^{\alpha -1}
e^{-t(\gamma +A)}f \, dt.
\end{equation}
\end{prop}

\begin{proof}
It is not hard to see that if $\gamma > 0$ then
$-(\gamma + A)$ is m-accretive on the Hilbert space
$H^s$ (with either periodic or Dirichlet boundary conditions). Let $u$ be an
element of the dense subspace $H^{s+2}\subset H^s$, then
\[
\langle (\gamma+A)u,u \rangle_s = \langle (\gamma+A)^{s+1} u,u \rangle_{L^2} =
\|u\|_{s+1}^2 \geq 0.
\]
Here we have used the inner product on $H^s$ as described in
Section \ref{sec:space}.

Standard results from the theory of elliptic boundary value problems also
give the existence of a solution $u \in H^{s+2}$ to the equation
$((1+\gamma)I +A)u = f$, subject to either of the two boundary conditions.
Thus $(1+\gamma)I +A$ is a surjection from $H^{s+2}$ to $H^s$ (this was
also discussed in Section \ref{sec:space}). This establishes the second
condition in (\ref{eqn:accretive}). Hence, by the Lumer-Phillips Theorem,
$-(\gamma+A)$ generates a contraction semigroup.

To show the integral in (\ref{eqn:fpow3}) converges in
$H^s$ for $\alpha \in (0,\infty)$ we establish the bound
$\|e^{-t(\gamma +A)}\|_{\mathcal{L}(H^s)}\leq e^{-\gamma t}$.
Since the above argument applies to any $\gamma >0$,
for any $n\in \mb{N}$ the operator $-(\gamma \frac{1}{n} +A)$ generates
a contraction semigroup on $H^s$. Thus
\[
\|e^{-t(\gamma +A)}\| =
\|e^{-t\gamma \frac{n-1}{n} }e^{ -t(\gamma \frac{1}{n}+A)}\| \leq
e^{-t\gamma \frac{n-1}{n} } \|e^{ -t(\gamma \frac{1}{n}+A)} \| \leq
e^{-t\gamma \frac{n-1}{n} }
\]
for $n$ arbitrarily large. Hence we can apply (\ref{eqn:fpowB}) to complete
the proof.
\end{proof}
We note that (\ref{eqn:fpow3}) is equivalent to (\ref{eqn:fpow2}) only for
$\alpha \in (0,1)$. The representation  (\ref{eqn:fpow3})
will be useful in proving a
comparison principle for the operator $(\gamma+A)^{-\beta}$ in Proposition
\ref{prop:Xcp}. We now show that $L$ generates a contraction semigroup on
$H^s$.

\begin{prop}
For each $s\geq 0$, the operator $L:=-(\gamma+A)^{1-\beta}$ is m-accretive
in $H^s$, and therefore generates the semigroup
$e^{tL} \in C([0,\infty),H^s)\cap C^1((0,\infty),H^s)$.
\end{prop}
This follows from a more general result in \cite{Kato60}, but the proof in this
case is short and straightforward, which is the following.
 \begin{proof}
We know $((1+\gamma)I + A)^{\alpha}$ maps
$H^{s+2\alpha}$ onto $H^s$ from the discussion in Section \ref{sec:space},
and for $u \in H^{s+\alpha}$ we have
\[
\langle (\gamma+A)^{\alpha}u,u \rangle_s = \langle (\gamma+A)^{s+\alpha} u,u
\rangle_{L^2} = \|u\|_{s+\alpha}^2 \geq 0.
\]
Hence, $-(\gamma+A)^{1-\beta}$ satisfies the hypotheses of the Lumer-Phillips
Theorem.
\end{proof}

We now use the results from (\ref{eqn:B1}), (\ref{eqn:B2}), and (\ref{eqn:B3})
 to show that $e^{tL}$ is a smoothing operator and to
establish operator bounds on $e^{tL}$ and $e^{tL}-I$. We say that $e^{tL}$ is
smoothing if it increases the regularity of a function as measured by
membership in different spaces. That is, an operator is smoothing
if it maps a space of functions into another space of smoother functions.

\begin{prop}
If $u \in L^2$, then for any $s>0$, $e^{tL} u \in H^s$. In particular,
for $n \in \mb{N}$, $\alpha \in (0,1)$, and $\lambda = 1-\beta$
we have the four bounds
\begin{equation} \label{eqn:smooth1}
\|e^{tL}\|_{\mathcal{L}(H^s,H^{s+2n\lambda})} \leq
\left(\frac{n}{\sqrt{2}t}\right)^n, \quad
\| (-L)^{\alpha} e^{tL} \|_{\mathcal{L}(H^s)} \leq C_{\alpha,T}
\left( \frac{\alpha}{t}\right)^{\alpha},
\end{equation}
\begin{equation} \label{eqn:smooth2}
\| e^{tL} \|_{\mathcal{L}(H^s,H^{s+2\alpha \lambda})} \leq
C_{\alpha,T} \left( \frac{\alpha}{t}\right)^{\alpha}, \quad
\| e^{tL}-I \|_{\mathcal{L}(H^{s+2\alpha \lambda},H^s)} \leq
C'_{\alpha,T}t^{\alpha}.
\end{equation}

\end{prop}
\begin{proof}
 We first note that because $L$ is self-adjoint and m-accretive on $H^s$ for
$s \geq 0$, estimate (\ref{eqn:B1}) gives
\begin{equation} \label{eqn:show}
\|L^n e^{tL}\|_{\mathcal{L}(H^s)} \leq \left( \frac{n}{\sqrt{2}t}\right)^n.
\end{equation}
Recall, as we established in Section \ref{sec:space}, that
$L:H^{s+2\lambda} \to H^s$ and that the inner product on $H^s$ is given by
$\langle u,v \rangle_s = \langle (\gamma + A)^s u,v \rangle_{L^2} =
\langle (-L)^{s/\lambda} u,v\rangle_{L^2}$, where we set $\lambda = 1-\beta$.
Note that $\lambda >0$. If $u \in H^s$, we compute
\begin{eqnarray*}
 \|e^{tL}u\|^2_{s+2n\lambda} &=& \langle e^{tL}u,e^{tL}u\rangle_{s+2n\lambda}
= \langle (-L)^{2n\lambda /\lambda}e^{tL}u, e^{tL} \rangle_s \\
&=& \langle (-L)^ne^{tL}u, (-L)^n e^{tL} \rangle_s =
(-1)^{2n}\|L^{n}e^{tL}u\|^2_s \leq
\left(\frac{n}{\sqrt{2}t}\right)^{2n} \|u\|^2_s.
\end{eqnarray*}
This, establishes the first bound in (\ref{eqn:smooth1}) and shows that
for $u \in L^2$, then $e^{tL} u \in H^s$ for any $s>0$. Hence
$e^{tL}$ is smoothing in the sense described above.

To apply (\ref{eqn:B3}) in the case $B=L$ and $H = H^s$, we have that
$Y = H^{s+2\alpha \lambda}$. Then estimates (\ref{eqn:B2}) and (\ref{eqn:B3})
yield the remaining three bounds in $(\ref{eqn:smooth1})$
 and $(\ref{eqn:smooth2})$.
\end{proof}
In fact, the bound on $\| e^{tL} \|_{\mathcal{L}(H^s,H^{s+2\alpha \lambda})}$
can be obtained from the bound on $\| (-L)^{\alpha} e^{tL} \|_{\mathcal{L}(H^s)}$
and a calculation similar to
that above for $\|e^{tL}u\|^2_{s+2n\lambda}$. In particular
\[
 \|e^{tL}u\|^2_{s+2\alpha \lambda}=\langle (-L)^{\alpha}e^{tL}u,
(-L)^{\alpha}e^{tL}u \rangle_s \leq C_{\alpha,T}^2
\left( \frac{\alpha}{t}\right)^{2\alpha}\|u\|_s^2.
\]

The final two results from this section provide comparison
principles for the operators $e^{tL}$ and $X$. They rely on
a comparison principle for the semigroup $e^{-t(\gamma +A)}$
and the integral formulas from Section \ref{sec:pow}.
The operator $\gamma+A$ is uniformly elliptic so the maximum principle for
parabolic equations applies to
\begin{equation} \label{eqn:ellmax}
\frac{\partial u}{\partial t} = -(\gamma+A)u, \quad u(0,x)=u_0(x).
\end{equation}
Thus, a solution to (\ref{eqn:ellmax}) on the interval $[0,T]$ obtains is
maximum on the boundary of $\Omega \times (0,T]$, but not at
$\Omega \times \{T\}$, (see \cite{Protter67}). Therefore, if $u \geq 0$,
then for $t>0$ we have $e^{-t(\gamma+A)}u\geq 0$. We now use this fact to
establish a comparison principle for $X$.
\begin{prop} \label{prop:Xcp}
 Let $\gamma > \sup_{x,y}|V_{22}(x,y)|$. Then $X$ satisfies a comparison
principle. That is, if $u\geq v$ a.e., then $X(u) = (\gamma+A)^{-\beta}
(\gamma u - V_2(x,u))\geq (\gamma+A)^{-\beta}(\gamma u - V_2(x,u)) = X(v)$ a.e.
\end{prop}
\begin{proof}
Since $\gamma > \sup_{x,y}\{|V_{22}(x,y)|\}$, then
$\gamma u - V_2(x,u)$ is increasing in $u$. Thus, $u\geq v$ implies
$\gamma u - V_2(x,u) \geq \gamma v - V_2(x,v)$. Then, as mentioned above, the
maximum principle for parabolic equations implies
\[
e^{-t(\gamma+A)}(\gamma u - V_2(x,u)) \geq e^{-t(\gamma+A)}( \gamma v - V_2(x,v)).
\]
Hence, for $t \geq 0$,
\[
t^{\beta -1}e^{-t(\gamma +A)}(\gamma u - V_2(x,u)) \geq t^{\beta -1}e^{-t(\gamma +A)}
( \gamma v - V_2(x,v)),
\]
thus the representation of $(\gamma+A)^{-\beta}$ in equation (\ref{eqn:fpow3})
implies that $(\gamma+A)^{-\beta}(\gamma u - V_2(x,u)) \geq
(\gamma+A)^{-\beta}(\gamma v - V_2(x,v))$, and we conclude that $X(u)\geq X(v)$.
\end{proof}

\begin{prop}\label{prop:Lcp}
If $u \geq v$ a.e. in $\Omega$, then $e^{tL}u \geq e^{tL}v$ in $\Omega$.
Moreover, we have the formula
\begin{equation} \label{eqn:Boch}
e^{tL} = e^{-t(\gamma+A)^{\lambda}} =  \int_0^{\infty} e^{-\tau(\gamma+A)}
\phi_{t,\lambda}(\tau) \, d\tau, \quad \forall t>0,
\end{equation}
with $\phi_{t,\lambda}(\tau)$ defined in (\ref{eqn:phi}).
\end{prop}
\begin{proof}
Equation (\ref{eqn:Boch}) is a valid application of (\ref{eqn:BochGen})
because $-(\gamma + A)$ is m-accretive.
If $u\geq v$ a.e. then for each $t>0$ and $\tau >0$, we have
\[
e^{-t(\gamma+A)}\phi_{t,\lambda}(\tau)u \geq e^{-t(\gamma+A)}\phi_{t,\lambda}
(\tau)v
\]
because $\phi_{t,\lambda} \geq 0$ and $e^{-t(\gamma+A)}$ satisfies
a comparison principle as explained above.
Integrating both sides of the inequality yields
$e^{tL}u \geq e^{tL}v$ by the subordination identity (\ref{eqn:Boch}).
The fact that $e^{tL}$ is smoothing,
in the sense of inequality (\ref{eqn:smooth2}), guarantees that
$e^{tL}u$ and $e^{tL}v$ are continuous from $\Omega$ to $\mb{R}$,
and therefore $e^{tL}u \geq e^{tL}v$ for all $x \in \Omega$.
\end{proof}

%=================================================
\subsection{$L^{\infty}$ bounds for $X$ and $e^{tL}$} \label{sec:bounds}

In preparation for the proof of existence and uniqueness of solutions to
(\ref{eqn:simple}), which will require $L^{\infty}$ estimates on $X$ and
$e^{tL}$, we will show that $X$ and $e^{tL}$ are, in fact, locally bounded
maps from
$L^{\infty}$ to itself. This is clear for $e^{tL}$ by the remark at the end of
Section \ref{sec:pow} because for $t>0$, $e^{tL}u \in H^s$ for
arbitrarily large $s$, and therefore it is in $L^{\infty}$ by the Sobolev
embedding theorem. However, we can use the comparison principles for
$X$ and $e^{tL}$ to provide explicit bounds.
\begin{prop} \label{prop:infty}
$X:L^{\infty} \to L^{\infty}$ is locally bounded with $\| X(u)\|_{L^{\infty}}
\leq \gamma^{\lambda}\|u\|_{L^{\infty}} + \gamma^{-\beta}\|V_2\|_{L^{\infty}}$.
Additionally, for each $t>0$, $e^{tL}: L^{\infty} \to L^{\infty}$ is a bounded
linear map and
$\|e^{tL}\|_{\mathcal{L}(L^{\infty})} \leq e^{-\gamma^{\lambda}t}$.
\end{prop}
\begin{proof}
Let $u \in L^{\infty}$, and set $\overline{u}=\|u\|_{L^{\infty}}$, so that
$u \leq \overline{u}$. Then the comparison principles established in
Propositions \ref{prop:Xcp} and \ref{prop:Lcp} imply that $X(u) \leq
X(\overline{u})$ and $e^{tL}u \leq e^{tL}\overline{u}$. Hence, if $X$ and
$e^{tL}$ are bounded on constant functions then they are bounded on
$L^{\infty}$.

Consider first $X$, and set $C_V = \gamma \overline{u} + \|V_2\|_{L^{\infty}}$.
We see that the boundedness of $V_2$ gives
$\gamma \overline{u} - V_2(x,\overline{u}) \leq C_V$. The integral
representation of $(\gamma+A)^{-\beta}$ in equation (\ref{eqn:fpow3}) and the
same reasoning as in the proof for Proposition \ref{prop:Xcp} imply
$X(\overline{u})=(\gamma+A)^{-\beta}(\gamma \overline{u} - V_2(x,\overline{u}))
\leq (\gamma+A)^{-\beta}C_V$. Thus, to establish bounds on $X$, we only need to
understand how
$(\gamma+A)^{-\beta}$ acts on constants.

A consequence of Proposition 10.3 from
\cite{Shubin78}, page 93, is that if $\psi$ is an eigenfunction of $\gamma+A$
with eigenvalue $\mu$, then $\psi$ is also an eigenfunction of
$(\gamma+A)^{-\beta}$ with eigenvalue $\mu^{-\beta}$. But  $(\gamma+A)C_V =
\gamma C_V$ because $A$ is a second-order differential operator. Hence
$(\gamma+A)^{-\beta}C_V = \gamma^{-\beta} C_V$ and
\[
X(u)\leq X(\overline{u})\leq \gamma^{-\beta} (\gamma \overline{u} +
\|V_2\|_{L^{\infty}})
=  \gamma^{\lambda}\|u\|_{L^{\infty}} + \gamma^{-\beta}\|V_2\|_{L^{\infty}},
\]
establishing the first claim in Proposition \ref{prop:infty}.

To bound $e^{tL}$ we examine how $L$ acts on constants. Using
$(\gamma+A)^{-\beta}\overline{u} = \gamma^{-\beta}\overline{u}$, we calculate
\[
L\overline{u}=-(\gamma+A)^{\lambda}\overline{u} = -(\gamma+A)(\gamma+A)^{-\beta}
\overline{u}= -(\gamma+A) \gamma^{-\beta}\overline{u} =
-\gamma^{1-\beta}\overline{u} = -\gamma^{\lambda}\overline{u}.
\]
 Now $e^{tL}$ is the semigroup generated by $L$, so we know
$e^{tL}\overline{u}$ solves
$\partial_t u = Lu$, $u(0) = \overline{u}$.  But $Le^{tL} = e^{tL}L$, hence
\[
\partial_t e^{tL}\overline{u} = e^{tL}L\overline{u}=e^{tL}
(-\gamma^{\lambda}\overline{u})= -\gamma^{\lambda} e^{tL}\overline{u}.
\]
Thus $e^{tL}\overline{u}$ solves $\partial_t u = -\gamma^{\lambda}u$,
$u(0)=\overline{u}$. Hence $e^{tL}\overline{u} = e^{-\gamma^{\lambda}t}
\overline{u}$, and we have $e^{tL}u \leq e^{tL}\overline{u} =
e^{-\gamma^{\lambda}t}\overline{u}$.
\end{proof}

%====================================================================
\section{Existence and uniqueness of solutions to
equation (\ref{eqn:simple})}\label{sec:eu}
The comparison result in Theorem \ref{thm:1} requires only the existence of the
flow generated by (\ref{eqn:simple}) for some short time $T >0$. However, the
motivation for studying this flow is to find critical points of the functional
(\ref{eqn:vp}), for which the flow must be defined on all of $(0,\infty)$.
In this section we establish the following
\begin{prop} \label{prop:ext}
If the potential $V\in C^{r+1}(\Omega \times \mb{R}, \mb{R})$, $r \geq 1$,
and $u_0 \in L^{\infty}(\Omega)$, then for every $\delta < 2$ there exists a
unique solution, $u(x,t) \in C((0,\infty),H^{r+\delta} \cap L^{\infty})$,
to (\ref{eqn:simple}) with $u(x,0)=u_0(x)$. If $r \geq 2$ then
$u(x,t) \in C^1((0,\infty),H^{r-2\lambda})$.
\end{prop}
We begin by showing existence and uniqueness of a mild solution in $L^{\infty}$,
and use the smoothing properties of the flow to obtain the desired regularity.
\begin{defn}
We say $u$ is a \emph{mild solution} of the equation
$\partial_t u = Lu + X(u)$, with $u(0,x) = u_0(x)$ if $u$ satisfies
\[
u(t,x) = e^{tL}u_0(x) + \int_0^t e^{(t-\tau)L}X(u(\tau,x)) \, d\tau.
\]
\end{defn}
We consider the map $\Psi$, on $C([0,T],L^{\infty})\, :
\,u(0,x)=u_0 \in L^{\infty}\}$ given by
\[
\Psi u (t) = e^{tL}u_0 + \int_0^t e^{(t-\tau)L}X(u(\tau)) \, d\tau.
\]
The mild solution, $u(x,t)$, of (\ref{eqn:simple}) will be the
fixed point of $\Psi$.
\begin{lem} \label{lem:psi}
For $T\in \mb{R}$, define $W_T=\{u \in C([0,T],L^{\infty})\, :
\,u(0,x)=u_0 \in L^{\infty}\}$. Then for small enough $T>0$,
$\Psi$ is a contraction on $W_T$. The size of $T$ is independent of $u_0$.
\end{lem}
\begin{proof}
The norm on $W$ will be defined as $\|u\|_{\infty,T} =
\sup_{\tau\in [0,T]}\|u(\tau)\|_{L^{\infty}}$. For any $t \in [0,T]$, it is
clear that $\|u\|_{\infty,t} \leq \|u\|_{\infty,T}$.
Let $u \in W_T$, then $\Psi u (0) = u_0$, and we have
\begin{eqnarray*}
\|\Psi u(t)\|_{L^{\infty}} &\leq& \|u_0\|_{L^{\infty}} + \int_0^t
e^{-\gamma^{\lambda} (t-\tau)} (\gamma^{\lambda} \| u(\tau)\|_{L^{\infty}} +
\gamma^{-\beta} \| V_2(x,u(x,\tau))\|_{L^{\infty}}) \, d\tau \\
&\leq&  \|u_0\|_{L^{\infty}} + (1-e^{-\gamma^{\lambda}t}) (\gamma^{\lambda}
\| u\|_{\infty,T} +\gamma^{-\beta} \| V_2\|_{L^{\infty}}) \leq C_T,
\end{eqnarray*}
where $C_T$ is finite for finite $T$. This follows directly from Proposition
\ref{prop:infty}.

Note that by the differentiability assumptions on $V$, we know that for any
$x \in \Omega$ and any $y_1,y_2 \in \mb{R}$, $|V_2(x,y_1)-V_2(x,y_2)| \leq
|V_2|_{C^1}|y_1-y_2|$. To see that $\Psi$ is a contraction, we compute
\begin{eqnarray*}
\| \Psi u(t) - \Psi v(t) \|_{L^{\infty}} & \leq & \int_0^t
\|e^{(t-\tau)L}(X(u)-X(v))\|_{L^{\infty}}\, d\tau \\
& \leq & \int_0^t e^{-\gamma^{\lambda} (t-\tau)} \| \gamma^{\lambda} (u - v)
+ \gamma^{-\beta} (V_2(x,v) -  V_2(x,u))\|_{L^{\infty}} \, d\tau \\
& \leq & (1-e^{-\gamma^{\lambda}t}) (\gamma^{\lambda} +
\gamma^{-\beta} |V_2|_{C^1})\|u - v\|_{\infty,T} = Ct\|u - v\|_{\infty,T},
\end{eqnarray*}
where $C=\gamma^{\lambda}(\gamma^{\lambda} + \gamma^{-\beta} |V_2|_{C^1})$
depends only on $\gamma$, $\beta$, and $V$.
Here we have used only for convenience the
fact  that $1-e^{-\gamma^{\lambda}t} \leq \gamma^{\lambda}t$ for $t\geq 0$.
Choosing $T = \frac{1}{2C}$ ensures that $\| \Psi u - \Psi v \|_{\infty,T} \leq
\frac{1}{2}\|u-v\|_{\infty,T}$, and thus $\Psi$ is a contraction on $W_T$.
This choice of $T$ depends depends only on $\gamma$, $\beta$,
and $V$. In particular, $T$ is independent of the initial condition
$u_0 \in L^{\infty}$.
\end{proof}

\begin{proof}[Proof of Proposition \ref{prop:ext}]
Lemma \ref{lem:psi} ensures that, for $T$ sufficiently small,
$\Psi$ has a unique fixed point $u \in C([0,T],L^{\infty})$ satisfying
\[
 u(t) = e^{tL}u_0 + \int_0^t e^{(t-\tau)L}X(u(\tau)) \, d\tau,
\]
establishing the existence and uniqueness of a mild solution to
(\ref{eqn:simple}) in $L^{\infty}$. Because $T$ was chosen independently of
$u_0$, we will also have existence on $[0,T]$ with initial condition $u(x,T)$,
which gives existence on the interval $[0,2T]$ for initial condition $u_0$.
This can be repeated indefinitely and we have $L^{\infty}$ existence on
$[0,\infty)$.

The operator $(\gamma+A)^{-\beta}$ is smoothing in the sense that it maps
$H^s$ into $H^{s+2\beta}$. Specifically, for $w \in H^s$,
\begin{eqnarray*}
\|(\gamma+A)^{-\beta}w\|_{s+2\beta}^2 &=& \langle (\gamma+A)^{-\beta}w,
(\gamma+A)^{-\beta}w \rangle_{s+2\beta} \\
&=& \langle (\gamma+A)^{-2\beta}w,w \rangle_{s+2\beta}
= \langle w,w \rangle_s = \|w\|_s^2,
\end{eqnarray*}
which gives $\|(\gamma+A)^{-\beta}\|_{\mathcal{L}(H^s,H^{s+2\beta})}=1$.
This implies that $X$ will have a smoothing property too. However,
this will be limited by the regularity of $V$.  Recall that $\Omega$
is bounded, so $L^{\infty} \subset L^2=H^0$, hence $(\gamma+A)^{-\beta}$
maps $L^{\infty}$ into $H^{2\beta}$, and for any $u \in L^{\infty}$, and
\[
\|X(u)\|_{2\beta} \leq \|\gamma u + V_2(x,u)\|_0 \, \leq \, \gamma \|u\|_{L^2}
+ \|V_2(x,u)\|_{L^2} \, \leq \,
|\Omega|^{1/2}(\gamma \|u\|_{L^{\infty}}+\|V_2\|_{L^{\infty}}).
\]
Thus $X(u(t))$ is bounded in $H^{2\beta}$ as long as $u(t)$ is bounded in
$L^{\infty}$.
Note that for any fixed $t\in [0,\infty)$ we know that
$\|u(\tau)\|_{\infty,t}$ is bounded, where the norm $\|\cdot \|_{\infty,t}$
is the same notation as in the proof of Lemma (\ref{lem:psi}). To see that
$u(t)$ is actually a solution in $H^{2\beta}$, we simply compute, for a fixed
$t>0$,
\begin{eqnarray*}
\|u(t)\|_{2\beta} &\leq & \|e^{tL}u_0\|_{2\beta} +
\int_0^t\|e^{(t-\tau)L}X(u(\tau))\|_{2\beta} \, d\tau \\
& \leq & \|e^{tL}u_0\|_{2\beta} +  |\Omega|^{1/2} \int_0^te^{-\gamma^{\lambda}
(t-\tau)}\left( \gamma \|u(\tau)\|_{L^{\infty}}+\|V_2\|_{L^{\infty}}\right)d\tau \\
& \leq & \|e^{tL}u_0\|_{2\beta} +  |\Omega|^{1/2} \frac{1}{\gamma^{\lambda}}
\left(1-e^{-\gamma^{\lambda}t}\right)
\left( \gamma \|u\|_{\infty,t}+\|V_2\|_{L^{\infty}} \right) < \infty.
\end{eqnarray*}

From here we wish to repeat this process to show that $u(t)$ is actually a
solution in $H^{4\beta}$. However, this will require a bound on the
composition $V(x,u(x,t))$ in the space $H^{2\beta}$. For this we will need to
employ the Moser estimates (\ref{eqn:me1}), but these estimates only apply for
$H^n$, $n\in \mb{N}$. This difficulty can be handled by splitting our analysis
into two cases, first for $\beta \in [1/2,1)$ and then for $\beta \in (0,1/2)$.
Recall that $e^{tL}$ is smoothing in the sense of estimates
(\ref{eqn:smooth1}) and (\ref{eqn:smooth2}),
so that if $u_0 \in L^{\infty} \subset L^2$ then for all $t>0$,
$e^{tL}u_0 \in H^s$ for any $s \geq 0$. Thus, the term $\|e^{tL}u_0\|_{n+1}$
that appears in the estimates below will not impede to our
regularity-building scheme.

Suppose $\beta \in [1/2,1)$, then the solution $u(t) \in L^{\infty}$ of
(\ref{eqn:simple}) is bounded in $H^{2\beta}$, but $2\beta \geq 1$, so
$\|u(t)\|_1 \leq \|u\|_{2\beta}$. Thus, $u(t)$ is an $H^1$-solution of
(\ref{eqn:simple}). To show that $u \in H^n$, we use induction on $n$, with
$n=1$ just established. Now assume $u \in H^n$ with $n\leq r$, then by
(\ref{eqn:me1}) we have $\|V_2(x,u)\|_n \leq c_V(1+\|u\|_n)$, where $c_V$
depends on $\|u\|_{\infty,t}$ and $V$. Because $2\beta \geq 1$ we have, for
any $w\in H^n$, $\|(\gamma+A)^{-\beta}w\|_{n+1} \leq \|w\|_{n+1-2\beta} \leq
\|w\|_n$. Using these facts, we compute
\begin{eqnarray*}
\|u(t)\|_{n+1}  & \leq & \|e^{tL}u_0\|_{n+1} + \int_0^t
\|e^{(t-\tau)L}X(u(\tau)) \|_{n+1} d\tau \\
& \leq & \|e^{tL}u_0\|_{n+1} + \int_0^t \|(\gamma+A)^{-\beta}
(\gamma u(\tau) - V_2(x,u(\tau)))\|_{n+1} d\tau \\
& \leq & \|e^{tL}u_0\|_{n+1} + \int_0^t  \|
\gamma u(\tau) - V_2(x,u(\tau))\|_n \, d\tau \\
& \leq &\|e^{tL}u_0\|_{n+1} + \int_0^t  \gamma \|u(\tau)\|_{n} + c_V
(1+\|u(\tau)\|_{n}) \, d\tau  \\
& \leq & \|e^{tL}u_0\|_{n+1} + Ct\sup_{\tau \in[0,t]}(1+\|u(\tau)\|_n)< \infty.
\end{eqnarray*}
Thus, $u(t) \in H^{n+1}$, for any $n \leq r$. Improved regularity to
$u(t) \in H^{r+\delta}$ for any $\delta <2$ will be shown after the case
$\beta <1/2$.

Between lines one and two of the above computation we have used the fact
that $e^{tL}$ is a contraction semigroup on $H^{n+1}$, i.e.
$\|e^{(t-\tau)L}\|_{\mathcal{L}(H^{n+1})} \leq 1$. We did not, however, take
advantage of the smoothing property of $e^{tL}$. This is something we cannot
afford to waste when $\beta \in (0,1/2)$.

Suppose $\beta \in (0,1/2)$, then $\lambda = 1-\beta \in (1/2,1)$. Let $u(t)$
be the $L^{\infty}$-solution of (\ref{eqn:simple}). The previous method of
using the smoothing properties of $(\gamma +A)^{-\beta}$ fails because the
smoothing factor of $2\beta$ is too small. Fortunately, this is precisely when
estimate (\ref{eqn:smooth2}) provides a large smoothing factor from $e^{tL}$.
Applying (\ref{eqn:smooth2}) with $\alpha = 1/2\lambda$ gives
\[
\|e^{tL}\|_{\mathcal{L}(H^s,H^{s+1})} \leq C_{\lambda,T}t^{-1/2\lambda}.
\]
The cost of the smoothing is the factor of $t^{-1/2\lambda}$.
However, $2\lambda >1$ so this is integrable on $[0,T]$, so we can proceed
with a similar argument to the case for $\beta \in [1/2,1)$. Just as we did
not use the smoothing properties of $e^{tL}$ in the previous case,
we do not need the smoothing properties
of $(\gamma +A)^{-\beta}$ in this case. To begin induction on $n$, we note that
$u(t)\in L^{\infty} \subset H^0$. If $u(t)\in H^n$ with $n \leq r$, then
\begin{eqnarray*}
\|u(t)\|_{n+1}  & \leq & \|e^{tL}u_0\|_{n+1} + \int_0^t
\|e^{(t-\tau)L}X(u(\tau)) \|_{n+1} d\tau \\
& \leq & \|e^{tL}u_0\|_{n+1} + \int_0^t C_{\lambda,T}t^{-1/2\lambda}
\|(\gamma u(\tau) - V_2(x,u(\tau))\|_{n} d\tau \\
& \leq &\|e^{tL}u_0\|_{n+1} + \tilde{C}\sup_{\tau \in[0,t]}
\left( 1+\|u(\tau)\|_{n}\right) t^{1-1/2\lambda} < \infty,
\end{eqnarray*}
where $\tilde{C}$ depends on $c_V$, $\gamma$, $\lambda$, and $T$. Thus,
$u(t)\in H^n$ for all $n \leq r+1$.

In either of the two cases, $\beta \geq 1/2$ and
$\beta < 1/2$, we have established that $u\in H^{r}$. From here we can
easily improve to $u\in H^{r+2\beta}$ because
\[
\|u(t)\|_{r+2\beta}  \leq \|e^{tL}u_0\|_{r+2\beta} + \int_0^t
\gamma \|u(\tau)\|_{r} + c_V(1+ \|u(\tau)\|_{r})  d\tau,
\]
regardless of the size of $\beta$. Now we can push a little further and
consider $\|u(t)\|_{r+2\beta+\epsilon}$ for some $\epsilon >0$. We can use
estimate (\ref{eqn:smooth2}) to obtain
\[
\|e^{(t-\tau)L}X(u(\tau))\|_{H^{r+2\beta+2\lambda\alpha}} \leq C^{\ast}
(t-\tau)^{-\alpha} (1+\|u(\tau)\|_{r}),
\]
provided $\alpha \in(0,1)$. The right side of the inequality is integrable
on $[0,t]$. Thus $u(t)\in H^{r+2\beta+2\lambda\alpha}$ for any $\alpha < 1$,
which gives $u(t)\in H^{r+\delta}$ for any $\delta < 2\beta + 2\lambda = 2\beta
+2(1-\beta) = 2$. Hence $u\in C([0,T],H^{r+\delta}\cap L^{\infty})$ for any
$\delta <2$.

To finish the proof of Proposition \ref{prop:ext}, we must show that the
derivative, $u_t$, of $u$ exists in $H^{r-2\lambda}$ so that $u$ is actually
$C^1$ and not just a mild solution. To see this, set $p=r-2\lambda$, and
consider
$R_h = \frac{1}{h} \|u(t+h) - u(t) - h u_t\|_{p}$, where
\[
u_t := Le^{tL}u_0 + \int_0^tLe^{(t-\tau)L}X(u(\tau)) \, d\tau + X(u(t)),
\]
Showing $R_h \to 0$ as $h\to 0$ will prove the desired result. We have
\begin{eqnarray*}
R_h &\leq & \|\frac{1}{h}(e^{hL} - I)e^{tL}u_0 - Le^{tL}u_0\|_{p} + \\
& + &\| \int_0^t \left(\frac{1}{h}(e^{hL}-I) - L\right)e^{(t-\tau)L}X(u(\tau))
\, d\tau \|_{p} + \\
& + & \| \frac{1}{h} \int_t^{t+h} e^{(t+h-\tau)L}X(u(\tau)) \,
d\tau - X(u(t))\|_{p}
\end{eqnarray*}
For brevity, we will refer to the three terms on the right-hand side of the
inequality as $I_1,I_2,$ and $I_3$. $I_1$ goes to zero with $h$
because $L$ is the generator of the $C_0$-semigroup $e^{tL}$ on $H^p$ and
the fact that $e^{tL}u_0 \in D(L)$ for all $t>0$. In $I_2$, we can replace
$(e^{hL}-I)$ by $\int_0^h Le^{\sigma L} d\sigma$ and we see that
\begin{eqnarray*}
I_2 &= &\| \int_0^t e^{(t-\tau)L}\left( \frac{1}{h}\int_0^h Le^{\sigma L}\,
d\sigma -L\right)X(u(\tau)) \, d\tau\|_{p} \\
& = & \| \int_0^t e^{(t-\tau)L} \frac{1}{h}\int_0^h \left(e^{\sigma L} -I \right)
LX(u(\tau)) \, d\sigma d\tau\|_{p} \\
& \leq & \int_0^t  \frac{1}{h} \int_0^h
\|e^{\sigma L} - I\|_{\mathcal{L}(H^{p+2\beta \lambda},H^{p})}
\|LX(u(\tau))\|_{p+2\beta \lambda} \, d\sigma d\tau\\
&\leq & \int_0^t  \frac{1}{h} \int_0^h C_2 \sigma^{\beta} \,d\sigma
\|LX(u(\tau))\|_{p+2\beta \lambda} \, d\tau \\
&\leq & C_2' h^{\beta}  \int_0^t \left(\|(u(\tau))\|_{r} +
c_V(\|u(\tau)\|_{\infty})(1+\|u(\tau)\|_{r}) \right) \, d\tau.
\end{eqnarray*}
We have used $\|e^{(t-\tau)L}\|_{\mathcal{L}(H^p)} \leq 1$ and,
between lines three and four, applied
(\ref{eqn:smooth2}) with $\alpha = \beta$.
The estimate in the final line of the calculation follows from the
bound $\|(\gamma+A)^{1-2\beta}w \|_s \leq \|w\|_{s+2(1-2\beta)}$,
which applies to the operator $LX(u)$. Hence,
 the fact that $p+2\beta \lambda + 2(1-2\beta) =
r-2\beta^2 < r$ allows the use of (\ref{eqn:me1}). We know that for
any fixed $t >0$, $\|u(\tau)\|_{r}$ and $\|u(\tau)\|_{\infty}$ are
bounded on $\tau \in [0,t]$, so $I_2$
can be made arbitrarily small for a suitable choice of $h$.

Without loss of generality, we assume that $h <1$ and we set\\
$M_1 = \max_{t\leq \tau \leq t+1} \|u(\tau)\|_{\infty}$ and
$M_2 = \max_{t\leq \tau \leq t+1} \|u(\tau)\|_{p}$. Let $\epsilon >0$, by
continuity there exists $h\in(0,1)$ such that $\|u(\tau)-u(t)\|_{p} < \epsilon$
for $|\tau-t|\leq h$.
Thus we have
\begin{eqnarray*}
I_3 & \leq & \frac{1}{h} \int_t^{t+h}\|e^{(t+h-\tau)L}
\left(X(u(\tau))-X(u(t))\right)\|_{p}
\, d\tau \\
&\leq & \frac{1}{h} \int_t^{t+h}\| \gamma(u(\tau)-u(t)) -
(V_2(x,u(\tau)) - V_2(x,u(t)))
\|_{p-2\beta} \, d\tau \\
& \leq & \frac{1}{h} \int_t^{t+h}  C_V(M_1)(1+2M_2)
\|u(\tau)-u(t)\|_{r-2} \\
& \leq & \frac{C_{V,M}}{h} \int_t^{t+h} \|u(\tau)-u(t)\|_{r-2}
 \leq  C_{V,M}' \epsilon \frac{1}{h}h \leq \epsilon
C_{V,M}',
\end{eqnarray*}
so that $I_3$ can be made arbitrarily small. We have again
used $\|e^{(t+h-\tau)L}\|_{\mathcal{L}(H^p)} \leq 1$ and that $p-2\beta = r-2$.
We have also used (\ref{eqn:me2}) between lines two and three above.
Thus we have shown $I_1$, $I_2$, $I_3 \to 0$ as $h \to 0$, establishing
$u \in C^1([0,\infty),H^{r-2\lambda})$, completing the proof of
Proposition \ref{prop:ext}.
\end{proof}

%\begin{rem}
%If $u_0 \in C^{r-1}$ then $u(t,\cdot) \in C^{r-1}$ as well. This can
%be seen from the representation of $u$ via Duhamel's formula
%(\ref{eqn:duhamel}).
%\end{rem}

%====================================================================
\section{Proof of Theorem \ref{thm:1}} \label{sec:comp}
We have established comparison principles for the semigroup $e^{tL}$ and the
operator $X$, as well as the existence of solutions to (\ref{eqn:simple}).
To emphasize the initial conditions, it will be convenient to write the
solutions of (\ref{eqn:simple}) as $u(x,t) = \Phi_tu_0$.
Hence, we aim to show that if $u_0$, $v_0 \in L^{\infty}$ and
$u_0 \geq v_0$ then $\Phi_t u_0 \geq \Phi_t v_0$ on a short time interval
$[0,T]$. This will follow from the iteration method below.

\subsection{Iteration method}
For $u \in L^{\infty}$ we define $F_t^0u = e^{tL}u$, and the $j$th iterate of
$u$ as
\[
 F_t^{j+1}u = e^{tL}u + \int_0^t e^{(t-\tau)L}X(F_{\tau}^ju) \, d\tau ,
\]
defined on some interval $[0,T]$. $F_t^{j+1}$ is well defined because
$X$ and $e^{tL}$ are both bounded maps from $L^{\infty}$ to itself, as shown in
Proposition \ref{prop:infty}.
\begin{prop} \label{prop:iter}
Let $T>0$. If $u \geq v$ then $F_t^ju \geq F_t^jv$ for all $t\in[0,T]$.
\end{prop}
\begin{proof}
Assume $u \geq v$. Then by Proposition \ref{prop:Lcp},
$F_t^0u = e^{tL}u \geq e^{tL}v = F_t^0v$. We assume that $F_t^ju \geq F_t^jv$
and proceed by induction on $j$. By Proposition \ref{prop:Xcp}, we have that
$X(F_t^ju) - X(F_t^jv) \geq 0$. Once again invoking Proposition \ref{prop:Lcp}
we have $e^{(t-s)L}[X(F_s^ju) - X(F_s^jv)] \geq 0$. Hence
\[
F_t^{j+1}u - F_t^{j+1}v = e^{tL}(u-v) + \int_0^t
e^{(t-s)L}\left[ X(F_s^ju) - X(F_s^jv) \right] \, ds \geq 0
\]
because the integrand is positive and $e^{tL}(u-v) \geq 0$.
\end{proof}
We need to show that this iteration converges to the solution in
Theorem \ref{thm:1}, so we now focus on a single initial condition, $u_0$. For
notational convenience we write $u^j(t)$ in place of $F_t^ju_0$, and
$u^0(t) = e^{tL}u_0$. Thus for each $j\in \mathbb{N}$ we have
\[
u^{j+1}(t) = e^{tL}u_0 + \int_0^t e^{(t-s)L}X(u^j(s)) \, ds.
\]
\begin{prop}
If $u_0\in L^{\infty}$ then there exists a $T>0$ such that for all $t\in(0,T]$,
$u^j(t) \in L^{\infty} \cap H^s$, for all $0\leq s \leq r+1$ and every
$j \in \mb{N}$.
\end{prop}
\begin{proof}
A consequence of Proposition \ref{prop:infty} is that $e^{tL}$ and $X$
are bounded on $L^{\infty}$. We use this to compute
\begin{eqnarray*}
\|u^{j+1}(t)\|_{\infty} & \leq & \|e^{tL}u_0\|_{\infty} + \int_0^t
\|e^{(t-s)L}X(u^j(s))\|_{\infty} \, ds \\
& \leq & \|u_0\|_{\infty} + T \max_{0 \leq s \leq T} \|u^j\|_{\infty} \\
& \leq & \|u_0\|_{\infty} + T(\|u_0\|_{\infty} + T \max_{0 \leq s \leq T}
\|u^{j-1}\|_{\infty} )\leq \ldots \\
& \leq &  \|u_0\|_{\infty} T + \|u_0\|_{\infty} T^2 +
\ldots +\|u_0\|_{\infty} T^j + T^{j+1}\|u^0\|_{\infty} \\
& \leq &  \|u_0\|_{\infty} \frac{1-T^{j+1}}{1-T} +T^{j+1}\|u_0\|_{\infty}
\leq \|u_0\|_{\infty} \frac{1}{1-T} = C_0
\end{eqnarray*}
Where, without loss of generality, we have assumed $T<1$. So we know that each
iterate $u^j$ is contained in the ball with radius $C_0$ in $L^{\infty}$ for
any $j$ and for all $t\in [0,T]$. The $L^{\infty}$ bounds allow the use of the
Moser estimates (\ref{eqn:me1}). For integer $k \leq r+1$ we have
\begin{eqnarray*}
\|u^{j+1}(t)\|_k & \leq & \|e^{tL}u_0\|_k + \int_0^t \|e^{(t-s)L}X(u^j(s)) \|_k
\, ds \\
& \leq & \|e^{tL}u_0\|_k + C_0T\max_{0\leq s \leq T}(1+ \|u^j(s)\|_k) \\
& \leq & \|e^{tL}u_0\|_k + C_0T(1+ \|u_0\|_k + C_0T\max_{0\leq s \leq T}
(1+ \|u^{j-1}(s)\|_k)) \leq \ldots \\
& \leq & \|e^{tL}u_0\|_k + C_0T(1+\|u_0\|_k) +  \ldots +(C_0T)^j(1+\|u_0\|_k)\\&&
+ (C_0T)^{j+1}(1+\|u_0\|_k) \\
& \leq & \sup_{0\leq t \leq T}\|e^{tL}u_0\|_k+ \frac{(1+\|u_0\|_k)}{1-C_0T} = C_k.
\end{eqnarray*} \end{proof}
We have assumed that $TC_0<1$ and $T<1$, which require only that
$T < \frac{1}{1+\|u_0\|_{\infty}}$.

\begin{prop} \label{prop:convHr}
$u^j$ converges in $C([0,T],H^r)$.
\end{prop}

\begin{proof}
\begin{eqnarray*}
\|u^{j+1}(t) - u^j(t)\|_r &\leq & \int_0^t\|e^{(t-s)L}
[X(u^j(s))-X(u^{j-1}(s))] \|_r \, ds \\
&\leq & t\max_{0 \leq s \leq t} \|X(u^j(s))-X(u^{j-1}(s))\|_r \\
&\leq & t \max_{0 \leq s \leq t} C|V|_rC_0(1\!+\!\|u^j(s)\|_r \!+\! \|u^{j-1}(s)\|_r)
\|u^j(s)\!-\!u^{j-1}(s)\|_r \\
&\leq & t \max_{0 \leq s \leq t} C^{\ast} \|u^j(s)-u^{j-1}(s)\|_r,
\end{eqnarray*}
where $C^{\ast}$ depends on $|V|_r$, $C_0$, and $C_r$. So we have
\[
\max_{0 \leq t \leq T} \|u^{j+1}(t) - u^j(t)\|_r  \leq C^{\ast}T
\max_{0 \leq t \leq T}\|u^j(s)-u^{j-1}(s)\|_r.
\]
We can set $T = \frac{1}{2C^{\ast}}$, which will ensure that the sequence
$u^j(t)$ is Cauchy in $H^r$. The argument for this is the following. First,
notice that if
$\|u^{j+1}-u^j\| \leq \frac{1}{2} \|u^j-u^{j-1}\|$, then
$\|u^{j+1}-u^j\| \leq (\frac{1}{2})^j \|u^1-u^0\|$. So for any $\epsilon >0$,
choose $N$ such that $(\frac{1}{2})^{N-1}\|u^1 - u^0\| < \epsilon$. Then we
have, for $m > n \geq N$,
\begin{eqnarray*}
\|u^m-u^n\| & \leq & \frac{1}{2} \|u^m-u^{m-1}\|+\frac{1}{2}
\|u^{m-1}-u^{m-2}\| + \ldots + \frac{1}{2} \|u^{n+1}-u^n\| \\
& \leq & \left[ \left( \frac{1}{2}\right)^{m-1}+ \left( \frac{1}{2}
\right)^{m-2} + \ldots +\left( \frac{1}{2} \right)^{n} \right] \|u^1-u^0\| \\
& \leq & \left( \frac{1}{2} \right)^{m-1}\left[ 1+2+2^2+\ldots +2^{m-n-1}
\right] \|u^1-u^0\| \\
& \leq &\left( \frac{1}{2} \right)^{m-1} 2^{m-n}\|u^1-u^0\| \leq
\left( \frac{1}{2} \right)^{n-1} \|u^1-u^0\| \leq \epsilon.
\end{eqnarray*}
\end{proof}
\begin{cor} \label{cor:1}
There exists a $T>0$ such that
$\Phi_t$ satisfies a comparison principle on the interval $[0,T]$.
\end{cor}
\begin{proof}
From Proposition \ref{prop:convHr} we have the existence of
$u^{\infty}\in C([0,T],H^r)$ such that $u^j \to u^{\infty}$ in $C([0,T],H^r)$.
This function $u^{\infty}(t,x) = \lim_{j\to \infty} F_t^ju_0(x)$ must satisfy
\[
u^{\infty}(t,x) = e^{tL}u_0 + \int_0^t e^{(t-s)L}X(u^{\infty}(s,x)) \, ds,
\]
and therefore, by Proposition \ref{prop:ext}, $u^{\infty}(t,x) = \Phi_t u_0$.
By Proposition
\ref{prop:iter}, if $u_0 \geq v_0$ then $F_t^ju_0 \geq F_t^jv_0$
and therefore $\Phi_t u_0 \geq \Phi_t v_0$ on $[0,T]$.
Thus we know $\Phi_t$ obeys a comparison principle on a small time interval
$[0,T]$.
\end{proof}
This establishes the comparison principle on a finite time interval $[0,T]$
and therefore concludes the proof of Theorem \ref{thm:1}.
To see that this comparison holds for all time $t >0$, we have the
following lemma.
\begin{lem} \label{lem:time}
If $u_0 \geq v_0$ for a.e. $x \in \Omega$, and there exist a
time $t_1 > 0$ for which $\Phi_{t_1} u_0(x) < \Phi_{t_1} v_0(x)$ on a set
of positive measure  then for every $t>0$, $\Phi_{t}u_0(x)<\Phi_{t}v_0(x)$
on a set of positive measure.
\end{lem}
\begin{proof}
Let $u_0(x) \geq v_0(x)$ for a.e. $x\in \Omega$.  Suppose there is a
first time $t_1$ such that
 $\Phi_{t_1}u_0(x)< \Phi_{t_1}v_0(x)$ on a set of positive
measure. Then on this set of positive measure
\begin{equation} \label{eqn:t1}
e^{t_1 L}(u_0(x)-v_0(x)) + \int_0^{t_1} e^{(t_1 - \tau)L}
(X(\Phi_{\tau}u_0(x)) - X(\Phi_{\tau}v_0(x)))\, d\tau <0.
\end{equation}
However, $e^{t_1 L}(u_0(x)-v_0(x)) \geq 0$ for all $x \in \Omega$ by Proposition
\ref{prop:Lcp}. We have by assumption that for all $\tau \in [0,t_1)$
$\Phi_{\tau}u(x) \geq \Phi_{\tau}v(x)$ for a.e. $x\in \Omega$. Hence, by
Proposition \ref{prop:Xcp} for all $\tau \in [0,t_1)$
$X(\Phi_{\tau}u(y)) - X(\Phi_{\tau}v(y)) \geq 0$ for a.e. $x\in \Omega$.
Again applying Proposition \ref{prop:Lcp} we have for any, $x \in \Omega$,
$e^{(t_1 - \tau)L}(X(\Phi_{\tau}u(y)) - X(\Phi_{\tau}v(y))) \geq 0$
for all $\tau \in [0,t_1)$. Thus for each $\tau \in [0,t_1)$
the integrand in (\ref{eqn:t1}) is a non-negative function on $\Omega$, and
therefore the integral (which is a function in $H^{r+1}$)
must be non-negative for a.e. $x \in \Omega$.
Therefore the left side of inequality (\ref{eqn:t1}) is the sum of the two
terms that are non-negative for a.e. $x$
and cannot be strictly negative on a set of positive measure.
\end{proof}

Combining Corollary \ref{cor:1} and Lemma \ref{lem:time} we
have that $\Phi_t$ satisfies a comparison principle on the interval $[0,T]$,
for $T>0$ and therefore $\Phi_t$ satisfies a comparison principle
on $[0,\infty)$.

%====================================================================
\section{Fractional elliptic equations} \label{sec:feqns}
We can extend the methods above to gradient descent equations for
energy functionals of the form
\[
S_{\alpha}(u) = \frac{1}{2}\langle u(x), A^{\alpha}(x) u(x) \rangle_{L^2} +
\int_{\Omega} V(x,u) dx,
\]
where $A$ is given by (\ref{eqn:A}) and $\alpha \in (0,1)$.
The Euler-Lagrange equation for $S_{\alpha}$ is
$A^{\alpha} u+ V_2(x,u) = 0$.

For a fixed $\alpha >0$, we use the inner product for the Sobolev
space $H^{\alpha}$ given by
\[
\langle u,v\rangle^{\ast}_{\alpha} = \langle (\gamma + A^{\alpha})u,
v \rangle_{L^2}
\]
we have the inner product on $H^{\alpha r}$ given by
\[
\langle u,v\rangle^{\ast}_{\alpha r} =
\langle (\gamma + A^{\alpha})^r u, v \rangle_{L^2}.
\]
Now consider $\alpha \in (0,1)$.
To calculate the Sobolev gradient of $S_{\alpha}$, we first note that
the derivative of $S_{\alpha}$ is
$DS_{\alpha}(u)\eta = \langle \eta, A^{\alpha} u+ V_2(x,u)\rangle_{L^2}$.
Thus, the Sobolev gradient of $S_{\alpha}$
in $H^{\alpha \beta}$, $\beta \in (0,1)$, is calculated as
\[ \begin{split}
DS_{\alpha}(u)\eta &= \langle \eta, A^{\alpha} u+ V_2(x,u)\rangle_{L^2} \\
& = \langle \eta, (\gamma + A^{\alpha})^{\beta}(\gamma + A^{\alpha})^{-\beta}
(A^{\alpha} u+ \gamma u - \gamma u+ V_2(x,u)) \rangle_{L^2}\\
&=  \langle \eta, (\gamma + A^{\alpha})^{-\beta}
(A^{\alpha} u+ \gamma u - \gamma u+ V_2(x,u)) \rangle^{\ast}_{\alpha \beta}\\
&=  \langle \eta, (\gamma + A^{\alpha})^{1-\beta}u -
(\gamma + A^{\alpha})^{-\beta}(\gamma u- V_2(x,u)) \rangle^{\ast}_{\alpha \beta}.
\end{split}\]
Hence, the gradient descent equation is
\begin{equation} \label{eqn:fdescent}
\partial_t u= -(\gamma + A^{\alpha})^{1-\beta}u +
(\gamma + A^{\alpha})^{-\beta}(\gamma u- V_2(x,u)).
\end{equation}
More concisely, we write $\partial_t u = \tilde{L}u + \tilde{X}(u)$ with
$\tilde{L}u : =-(\gamma + A^{\alpha})^{1-\beta}u $ and
$\tilde{X}(u):=(\gamma + A^{\alpha})^{-\beta}(\gamma u- V_2(x,u))$.

Recall that the maximum principle for parabolic equations ensures that
the semigroup $e^{-tA}$ satisfies a comparison principle.
The Bochner subordination identity (\ref{eqn:BochGen}) allows us to write
the semigroup generated by  $-A^{\alpha}$ as
\begin{equation}
e^{-tA^{\alpha}} =  \int_0^{\infty} e^{-\tau A}
\phi_{t,\alpha}(\tau) \, d\tau, \quad t>0.
\end{equation}
Thus, by the same argument as in the proof of Proposition \ref{prop:Lcp},
the comparison principle for $e^{tA}$ guarantees that
$e^{-tA^{\alpha}}$  satisfies a comparison principle as well.

Analogously, $-(\gamma  + A^{\alpha})$ generates the semigroup
$e^{-t(\gamma +A^{\alpha})}$, which also satisfies a comparison principle.
Just as in (\ref{eqn:fpow3}) we can define the real powers as
\[
(\gamma +A^{\alpha})^{-\beta} f = \frac{1}{\Gamma (\beta)}
\int_0^{\infty}t^{\beta -1} e^{-t(\gamma +A^{\alpha})}f \, dt,
\]
and positive powers are again the inverses of negative powers. Thus, because
we have a comparison principle for $e^{-t(\gamma +A^{\alpha})}$ we have that
$(\gamma +A^{\alpha})^{-\beta}$ also satisfies a comparison principle.
An argument as in the proof of Proposition \ref{prop:Xcp} shows that if
$\gamma > \sup_{x,y} |V_{22}(x,y)|$ then $\tilde{X}$ will satisfy a comparison
principle.

Employing the subordination identity once again, we have
\[
e^{t\tilde{L}}=e^{-t(\gamma +A^{\alpha})^{1-\beta}} =
\int_0^{\infty} e^{-\tau(\gamma + A^{\alpha})}
\phi_{t,1-\beta}(\tau) \, d\tau, \quad t>0.
\]
Hence, we see $e^{t\tilde{L}}$ satisfies a comparison principle
because $\phi_{t,1-\beta}(\tau)>0$ for all $\tau >0$ and
 $e^{-\tau(\gamma + A^{\alpha})} $ satisfies comparison principle, just as in the
proof of Proposition \ref{prop:Lcp}.

Smoothing estimates like (\ref{eqn:smooth2}) and (\ref{eqn:smooth2})
follow for $e^{t\tilde{L}}$ just as they did for $e^{tL}$ in Section
\ref{sec:pow}.
The existence of solutions for $\partial_t u = \tilde{L}u + \tilde{X}(u)$
follows from these smoothing estimates and the arguments
from Section \ref{sec:eu}.

Finally, the proof for a comparison principle for the flow defined by
$\partial_t u = \tilde{L}u + \tilde{X}(u)$ follows from the Duhamel formula
\[
u(t,x) = e^{-t\tilde{L}}u_0(x) +
\int_0^t e^{-(t-\tau)\tilde{L}}\tilde{X}(u(\tau,x)) \, d\tau,
\]
and the iteration argument from Section \ref{sec:comp}.

%====================================================================
\subsection{Constant coefficients}
In the case of periodic boundary conditions and if the matrix $a(x)$ is
constant (i.e. independent of $x$), the situation is
simplified because we can write down concrete formulae for the
operators, and we can use classical Fourier analysis in place of some of the
abstract semigroup theory. For instance, we could have avoided
the theory of fractional powers of operators and instead used the definition
$(\gamma + A)^{\alpha}u =
((\gamma + 4\pi^2 \xi^T a\xi)^{\alpha}\hat{u}(\xi))^{\vee}$

One can also write down the semigroup $e^{-t(\gamma+A)}=e^{-t(\gamma-\Delta)}=
e^{-\gamma t}e^{t\Delta}$. The operator $e^{t\Delta}$ is convolution with the
heat kernel, written as
\begin{equation} \label{eqn:heat}
\left(e^{t\Delta}u\right)(x) =\frac{1}{(4\pi t)^{d/2}}\int_0^1
\sum_{k\in \mb{Z}^d}e^{-|x-y+k|^2/4t} u(y) \, dy.
\end{equation}
Combining this with equation (\ref{eqn:Boch}) yields the comparison principle
for $e^{tL}$ immediately.

%====================================================================
\section{Application to Aubry-Mather theory for PDEs} \label{sec:app}
We can use the comparison principle from Theorem \ref{thm:1} and the results
from \cite{delaLlaveVal08} to further develop the Aubry-Mather theory
for PDEs to the case
of a general elliptic problem of the form (\ref{eqn:ell}). For
this section, we restrict
discussion to periodic boundary conditions. We require the potential
function $V$ and the matrix coefficient functions $a^{ij}(x)$ to be periodic
over the integers.  That is,
\[
V(x+e,y +l) = V(x,y) \quad \forall (e,l) \in \mb{Z}^d \times \mb{Z},
\quad \forall (x,y) \in \mb{R}^d\times \mb{R}
\]
\[
a^{ij}(x+e) = a^{ij}(x) \quad \forall e \in \mb{Z}^d,
\quad \forall x \in \mb{R}^d
\]
which we write as $V:\mb{T}^d\times \mb{T} \to \mb{R}$ and $a^{ij}:\mb{T}^d \to
\mb{R}$. An important class of functions is given by the following
\begin{defn}
 A function $u: \mathbb{R}^d \to \mathbb{R}$ is said to have the
\emph{Birkhoff property} (or $u$ is a \emph{Birkhoff function}) if for any
fixed $e \in \mb{Z}^d$, and $ l \in \mb{Z}$, $u(x+e) - (u(x) + l)$ does not
change sign with $x$. That is, $u(x+e) - (u(x) + l)$ is either $\geq 0$ or
$\leq 0$ depending on the choices of $e$ and $l$, but not $x$. Any such
function can be seen as a surface in $\mb{T}^{d+1}$ without any self-crossings.
\end{defn}
For a fixed $\omega \in \mb{R}^d$, we define
$B_{\omega} = \{u : \mb{R}^d\to \mb{R} | u \mbox{ is Birkoff}, u(x)-\omega
\cdot x \in L^{\infty}(\mb{R}^d)\}$. We note that $u \in B_{\omega}$
if and only if $u$ is a Birkhoff function and if for any
$e\in \mb{Z}^d $ and $l \in \mb{Z}$, either
$u(x+e)-u(x) - l \leq 0$ or $\geq 0$
according to whether $\omega \cdot e -l \leq 0$ or $\geq 0$
(see \cite{Moser86}). The vector $\omega$ is referred to as the
\emph{frequency} or the \emph{rotation vector} of the function $u$,
and is a natural generalization of the one dimensional notion of rotation
number. We have the following result.
\begin{thm} \label{thm:2}
 Let $V \in C^2(\mb{T}^d \times \mb{T}, \mb{R})$, and let $A$ be a
self-adjoint, uniformly elliptic operator given by $Au=-\mbox{div} (a(x)
\nabla u)$ with coefficients $a^{ij} \in C^{\infty}(\mb{T}^d,\mb{R})$.
Then for any $\omega \in \mb{R}^d$, there exists a solution $u \in B_{\omega}$
to equation (\ref{eqn:ell}). That is  $A u + V_2(x,u) = 0$,
$u(x) - \omega \cdot x \in L^{\infty}(\mb{R}^d)$, and $u$ is Birkhoff.
\end{thm}

The method of proof from \cite{delaLlaveVal08} is to first show the result
holds for rational frequencies (i.e. for any $\omega_N \in \frac{1}{N}\mb{Z}^d$
with $N \in \mb{N}$). Then we obtain solutions for arbitrary
$\omega \in \mb{R}^d$ as limits of solutions with rational frequencies.
Passing to the limit requires an oscillation lemma of De Giorgi-Moser
type, as well as classical $C^{\epsilon}$ elliptic estimates.

For minimal solutions, the oscillation lemma was given by Moser in
Theorem 2.2 of \cite{Moser86}. This gives a bound on the supremum of
the solution $u_N$ by the norm of its associated frequency
$\omega_N \in \frac{1}{N}\mb{Z}^d$, and independent of $N$.
We provide the arguments for the oscillation lemma in Section
\ref{sec:irrat}.

To have $\omega_N \to \omega$, we must have $N \to \infty$, and
the size of the fundamental domain, $N\mb{T}^d$, becomes unbounded. If
the bound on $u_N$ depended on the size of the domain, then the
$C^{\epsilon}$-estimates on $\nabla u_N$ would degenerate. Instead, because
the bounds on the $u_N$ are in terms of $\omega_N$,
the $C^{\epsilon}$-estimates on $\nabla u_N$ are uniform in $N$.
Hence we can conclude the convergence of $u_N$ to a continuous function $u$
with associated frequency $\omega$.

Note that the regularity assumption on $V$ is weaker than in Theorem
\ref{thm:1} because we can choose a smooth initial condition
(i.e. $u_0(x) = \omega \cdot x$).

\subsection{Rational frequencies}
\begin{lem}
If $\omega \in \frac{1}{N}\mb{Z}^d$ and
$u(x,t)$ solves (\ref{eqn:general}) on $N\mb{T}^d$ with initial condition
$u_0(x)=\omega \cdot x$, then there exists a sequence $t_n \to \infty$ as
$n\to \infty$ such that $u(x,t_n) \to u^{\ast}_{\omega}$ in $L^2$ and
$u^{\ast}_{\omega}$ solves (\ref{eqn:ell}).
\end{lem}
\begin{proof}
Suppose $u(x,t)$ solves (\ref{eqn:general}) on $[0,N]^d$ with periodic boundary
conditions and initial condition $u_o = \omega \cdot x$. In Section
\ref{sec:space} we showed that
$DS_N(u)\eta= \langle (\gamma + A)^{1-\beta}u - (\gamma +A)^{-\beta}
\left(\gamma u - V_2(x,u) \right), \eta \rangle_{\beta}$. This together with
the assumption that $u(x,t)$ solves (\ref{eqn:general}) allows us
to show that $S_N(u(x,t))$ is decreasing in $t$. More precisely,
\begin{eqnarray*}
\frac{d}{dt}S_N(u(t))&=& DS_N(u(t))\partial_tu \\
&=& \langle (\gamma + A)^{1-\beta}u - (\gamma +A)^{-\beta}\left(\gamma u -
V_2(x,u) \right), \partial_tu \rangle_{H^{\beta}(N\mb{T}^d)}\\
&=& - \|(\gamma + A)^{1-\beta}u - (\gamma +A)^{-\beta}\left(\gamma u -
V_2(x,u) \right) \|_{H^{\beta}(N\mb{T}^d)}^2 \leq 0.
\end{eqnarray*}
Recalling that $\Lambda_1 \leq a(x) \leq \Lambda_2$, we have
$S_N(u(t)) \leq S_N(u_0)\leq \frac{N^d}{2}\Lambda_2 |\omega|^2 +
N^d\|V\|_{L^{\infty}}$ for all $t >0$. We can conclude that
\begin{equation}\begin{split}
\frac{\Lambda_1}{2}\int_{[0,N]^d} |\nabla u(x,t)|^2 dx &\leq
\int_{[0,N]^d} \frac{1}{2}a(x)\nabla u(x,t)\cdot \nabla u(x,t) dx \\
&\leq \frac{N^d}{2}\Lambda_2 |\omega|^2 + 2N^d\|V\|_{L^{\infty}([0,N]^d)}.
\end{split}\end{equation}
Thus, $\|\nabla u(t)\|_{L^2}$ is bounded uniformly in $t$.

Because $S_N(u)$
is bounded below and $\frac{d}{dt}S_N(u(t_n))\leq 0$, there is a sequence
$t_n \to \infty$ such that $\frac{d}{dt}S_N(u(t_n)) \to 0$ as $n \to \infty$.
The periodicity of $V$ ensures that for any sequence of integers, $\{k_n \}$,
we have $S_N(u(t_n)+k_n) = S_N(u(t_n))$. By selecting
$k_n = - \lfloor N^{-d}\int_{[0,N]^d}u(x,t_n)dx \rfloor $, we can replace
the sequence $u(t_n)$ by $u(t_n) +k_n$ without affecting $S_N$. Therefore
we assume, without loss of generality, that for each $n$, the average
of $u(x,t_n)$ lies in the range $[0,1]$. The gradients $\nabla u(x,t_n)$
are uniformly bounded in $L^2$, so by Poincar\'e's inequality we
have that $\{u(t_n)\}$ is precompact, and there is a subsequence $u(t_{n_k})$
that converges weakly in $H^1$ and strongly in $L^2$. Denote the limit by
$u^{\ast}_{\omega}$. $V_2$ is Lipschitz, so
$\|V_2(\cdot ,u(t_{n_k}))-V_2(\cdot , u^{\ast}_{\omega})\|_{L^2} \to 0$
as $k \to \infty$. $A$ is an $L^2$-closed operator, hence
$Au^{\ast}_{\omega} +V_2(x,u^{\ast}_{\omega})=0$.
\end{proof}
\begin{lem}
If $u(x,t)$ solves (\ref{eqn:general}) with initial condition
$u(x,0)=u_0(x)\in B_{\omega}$, then $u(t)\in B_{\omega}$ for all $t>0$.
\end{lem}
\begin{proof}
As before, we denote the solution $u$ with initial condition $u_0$ of
(\ref{eqn:general}) as $u(x,t)=\Phi_tu_0(x)$. We can conclude from
Theorem \ref{thm:1} that if $u_0 \leq v_0$ then $\Phi_tu_0(x)\leq \Phi_tv_0(x)$.
For convenience we let $\mc{C}_k$ and $\mc{R}_l$ denote the family of
operators
\begin{displaymath}
\mc{C}_ku(x)=u(x+k) \quad
\mbox{and} \quad \mc{R}_lu(x)=u(x)+l
\end{displaymath}
for each $k\in \mb{Z}^d$ and $l\in \mb{Z}$.

Let $l\in \mb{Z}$ and define $Y_l(x,t)=\mc{R}_l\Phi_tu_0(x)=u(x,t) + l$.
Recall the abbreviated notation
\begin{displaymath}
Lu=-(\gamma + A)^{1-\beta}u \quad \mbox{and} \quad
X(u)= (\gamma +A)^{-\beta}(\gamma u - V_2(x,u)),
\end{displaymath}
which allows us to write
(\ref{eqn:general}) as $\partial_tu = Lu + X(u)$, or
$\partial_t\Phi_tu_0 = L\Phi_tu_0 + X(\Phi_tu_0)$.

Now, $V_2$ is periodic over the integers, so
$V_2(x,u)=V_2(x,\mc{R}_lu)$. Also, $(\gamma +A)^{\alpha}c = \gamma^{\alpha}c$
for any constant $c\in \mb{R}$. Hence
\[ \begin{split}
X(Y_l)&=X(\mc{R}_l\Phi_tu_0)=X(u(x,t)+l)=
(\gamma +A)^{-\beta}(\gamma u+\gamma l - V_2(x,u+l))\\
&=\gamma^{1-\beta}l+(\gamma +A)^{-\beta}(\gamma u - V_2(x,u))=
\gamma^{1-\beta}l+ X(\Phi_tu_0).
\end{split}\]
Similarly,
\[ \begin{split}
LY_l&=L(\mc{R}_l\Phi_tu_0)=L(u(x,t)+l)=-(\gamma + A)^{1-\beta}(u(x,t)+l)\\
&=-\gamma^{1-\beta}l-(\gamma + A)^{1-\beta}u(x,t)=-\gamma^{1-\beta}l +
L(\Phi_tu_0).
\end{split}\]

Therefore $L\Phi_tu_0 + X(\Phi_tu_0)= LY_l  +X(Y_l)$ and we have
$\partial_tY_l=\partial_t\Phi_tu_0 = L\Phi_tu_0 + X(\Phi_tu_0)=
LY_l  +X(Y_l)$.
So $Y_l$ solves (\ref{eqn:general}) with initial condition $Y_l(x,0)=u_0(x)+l$.
But $\partial_t(\Phi_t\mc{R}_lu_0) =L(\Phi_t\mc{R}_lu_0)+X(\Phi_t\mc{R}_lu_0)$
with the same initial condition $\Phi_0\mc{R}_lu_0=u_0+l$. Thus,
\begin{equation} \label{eqn:R}
\Phi_t\mc{R}_lu_0=\mc{R}_l\Phi_tu_0
\end{equation}
by the uniqueness of solutions to
(\ref{eqn:general}) as shown in Proposition \ref{prop:ext}.

Now define $Z_k(x,t)=\mc{C}_k\Phi_tu_0 = u(x+k,t)$. The periodicity of $V_2$
ensures that $\mc{C}_kV_2(x,u(x,t))=V_2(x+k,u(x+k,t))=V_2(x,u(x+k,t)$ so that
$X(\mc{C}_ku(x,t))=\mc{C}_kX(u(x,t))$. Clearly
$L(\mc{C}_ku(x,t))=\mc{C}_kLu(x,t)$, so
\[\begin{split}
\partial_tZ_k &= \partial_tu(x+k,t)
 = \mc{C}_kLu(x,t)+\mc{C}_kX(u(x,t))\\
&=L\mc{C}_ku(x,t)+X(\mc{C}_ku(x,t))=
LZ_k + X(Z_k).
\end{split}\]
So $Z_k$ solves (\ref{eqn:general}) with initial condition
$Z_l(x,0)=u_0(x+k)$. But
$\partial_t(\Phi_t\mc{C}_ku_0) =L(\Phi_t\mc{C}_ku_0)+X(\Phi_t\mc{C}_ku_0)$
with the same initial condition $\Phi_0\mc{C}_ku_0=u_0(x+k)$. Thus,
\begin{equation} \label{eqn:C}
\Phi_t\mc{C}_ku_0=\mc{C}_k\Phi_tu_0
\end{equation}
by the uniqueness of solutions to
(\ref{eqn:general}) as shown in Proposition \ref{prop:ext}.

Now suppose $u_0 \in \mc{B}_{\omega}$ so that $\mc{C}_k\mc{R}_lu_0 \leq 0$
or $\geq 0$ according to whether $\omega \cdot k +l \leq 0$ or
$\omega \cdot k +l\geq 0$.Then the comparison principle from Theorem
\ref{thm:1} yields $\Phi_t\mc{C}_k\mc{R}_lu_0 \leq 0$ or $\geq 0$
according to whether $\omega \cdot k +l \leq 0$ or $\omega \cdot k +l\geq 0$.
But equations (\ref{eqn:R}) and (\ref{eqn:C}) together give
$\Phi_t\mc{C}_k\mc{R}_l = \mc{C}_k\mc{R}_l\Phi_t$, so
$\mc{C}_k\mc{R}_l\Phi_tu_0 \leq 0$ or $\geq 0$
according to whether $\omega \cdot k +l \leq 0$ or $\omega \cdot k +l\geq 0$,
hence $\Phi_tu_0 \in \mc{B}_{\omega}$.
\end{proof}

The set $\mc{B}_{\omega}$ is closed  under $L^2$-limits, so
$u^{\ast}_{\omega} \in \mc{B}_{\omega}$. This establishes Theorem \ref{thm:2}
in the case $\omega \in \frac{1}{N}\mb{Z}^d$.

\subsection{Irrational frequencies} \label{sec:irrat}
For the case of irrational frequency, let
$\omega \in \mathbb{R}^d\setminus \mathbb{Q}^d$ and let $(\omega_n )$
be a sequence such that $\omega_n \in \frac{1}{n}\mathbb{Z}^d$ for each $n$
and $\omega_n \to \omega$. Let $\mbox{div}(a(x)\nabla u_n) = V_2(x,u_n)$,
where $u_n \in \omega_n \cdot x + H^1(N\mathbb{T}^d)$, $\lambda < a(x) <
\Lambda$, and let  $B_R$ denote a ball of radius $R$ centered at some point
in the domain $\Omega$ such that the concentric ball $B_{4R}$ is also in the
domain $\Omega$. Then Theorem 8.22 of \cite{GilbargTrudinger01}
gives for each $n$ that
\[
\mbox{osc}_{B_R} (u_n) \leq (1-C^{-1}) \mbox{osc}_{B_{4R}}(u_n) + k(R)
\]
where $k=\lambda^{-1}R^{2(1-d/q)}\|V_2(x,u_n)\|_{L^{q/2}(\Omega)}$ for $ q>d$,
which is bounded by the volume of $\Omega$ times $\|V_2\|_{L^{\infty}}$, and
$C=C(d,\Lambda/\lambda,q)>1$ (pages 200-201 of \cite{GilbargTrudinger01}).
Following Moser's methods from the proof of Theorem 2.2 of \cite{Moser86},
we can obtain a local result on the cube
$Q = \{ x\in \mathbb{R}^d : |x_j| \leq 1/2 \}$, and then consider translations.

Using a variant of the above inequality we can say
\[
\mbox{osc}_{Q} (u_n) \leq \theta \mbox{osc}_{4Q}(u_n) + c_1
\]
where $\theta \in (0,1)$ depends on $\lambda$, $\Lambda$, $d$, since we can
fix $q=d+1$, and $c_1$ depends on $\lambda$, $d$, and
$\|V_2(x,u_n)\|_{L^{q/2}(\Omega)}$. Here we will take $\Omega = 5Q = \{|x_j|
\leq 5/2\}$ so that $\|V_2(x,u_n)\|_{L^{q/2}(\Omega)}
\leq 5^d\|(V_2)^{(d+1)/2}\|_{L^{\infty}}$, independent of $n$.

Since $u_n(x+k)-u_n(x)-l$ has the same sign as
$\omega_n \cdot k -l$, for any $x\in \mathbb{R}^d$ and any
$k\in \mathbb{Z}^d$, $l \in \mathbb{Z}$, we can conclude
$|u_n(x+k)-u_n(x)-\omega_n\cdot k | \leq 1$ for any such $x$ and $k$. For an
arbitrary $y \in \mathbb{R}^d$ select $k \in \mathbb{Z}^d$ such that
$y-k \in Q$, then we have
\[ \begin{split}
|u_n(x+y)-u_n(x)-y\cdot \omega_n| &\leq |u_n(x+y)-u_n(x+k)| \\
&\ \ +|u_n(x+k)-u_n(x)- \omega_n\cdot k| + | \omega_n \cdot k- \omega_n \cdot y| \\
& \leq \mbox{osc}_{x+k+Q}(u_n) + 1 + |\omega_n|
\end{split}\]
Then an argument parallel to the one on page 240 of \cite{Moser86}
shows that $\mbox{osc}_{4Q}(u_n) \leq
\mbox{osc}_Q(u_n) + 2\sum_{j=1}^d(1+2|(\omega_n)_j|)$ so that with the result
from \cite{GilbargTrudinger01} we have
$\mbox{osc}_Q(u_n) \leq c_2\sqrt{1+|\omega_n|^2}$. This
estimate holds in translated cubes $x+Q$, so if we assume that $|\omega_n|
\leq c |\omega|$, then the $|u_n(x+y)-u_n(x)-\omega_n \cdot y|$ are bounded in
$\mathbb{R}^d$ by $c_3\sqrt{1+|\omega|^2}$, independently of $n$.

This bound is crucial, because it allows the application of Theorem 5.2,
page 277 of \cite{LadyUral68}, which gives
$u_n \in \omega_n \cdot x+ C^{1,\epsilon}$ (Assuming $V\in C^{2,\epsilon}$), and
bounds
\[
|\nabla u_n| \leq \gamma_1
\]
with $\gamma_1$ depending on the ellipticity constants and $\omega_n$ (in fact,
 monotone in $\omega_n$), so bounded for our sequence of $(\omega_n)$. So,
convergence of $\omega_n \to \omega$ and the Arzel\`a-Ascoli Theorem imply the
existence of a subsequence $u_{n_k} \to u^{\ast}_{\omega}$ in $C^0_{loc}$.
Then, $\mbox{div}(a(x) \nabla )$ is closed under $C^0_{loc}$ limits, and we have
\[
Au^{\ast}_{\omega}+V_2(x,u^{\ast}_{\omega}) =
-\mbox{div}(a(x) \nabla u^{\ast}_{\omega})+V_2(x,u^{\ast}_{\omega})=0.
\]
This completes the proof of Theorem \ref{thm:2}. \qed
\begin{rem}
If we replace the operator $A$ from Theorem \ref{thm:2} with $A^{\alpha}$ for
$\alpha >0$, then the same method of proof for Theorem \ref{thm:2} will allow
us to find solutions of $A^{\alpha}u + V_2(x,u)=0$, with $u \in \mc{B}_{\omega}$
provided that $\omega$ is rational. However, to find solutions for irrational
$\omega$, we require a oscillation lemma for the critical points with
rational frequencies. With this, one may be able to find solutions with
irrational frequencies using methods from \cite{CabreSola05}.
The work of  \cite{CaffarelliSylvestre09} establishes a Harnack inequality in
the autonomous case, and can perhaps be extended to our case as well.

For this method to work, one would need both the oscillation lemma of
De Giorgi-Moser type as well as $C^{\epsilon}$ estimates for the fractional
power of an elliptic operator.
\end{rem}

\section*{Acknowledgements}
TB and RL were supported by NSF grants and by 
Texas Coordinating Board NHARP 0223. EV was supported by
GNAMPA and FIRB. We would like to thank J.W. Neuberger for
several helpful conversations.

\nocite{Taylor96} \nocite{Taylor97} \nocite{Miklavicic98}
\nocite{Stein70} \nocite{DeGiorgi57} \nocite{Haase06} \nocite{Martinez01}

%\bibliographystyle{foo}%{alpha}
%\bibliography{comparison}

\def\cprime{$'$}

\medskip

Received January 2010; revised June 2010.

\medskip

\end{document}